\documentclass[11pt, final]{article}
\usepackage[top=0.7in, bottom=0.7in, left=1.2in, right=1.2in]{geometry}
\usepackage{xcolor}
\usepackage[english]{babel}
\usepackage{amsmath,amssymb,amsthm,mathdots}
\usepackage{graphicx}
\usepackage[colorinlistoftodos]{todonotes}

\usepackage{xcolor}

\newcommand{\mb}[1]{\mathbb{#1}}
\newcommand{\mc}[1]{\mathcal{#1}}

\newcommand{\norm}[2]{\|#1\|_{#2}}

\newcommand{\normtwo}[1]{\left\|#1\right\|_{2}}

\newcommand{\B}[1]{\boldsymbol{#1}•}

\newcommand{\bmat}[1]{\begin{bmatrix}#1\end{bmatrix}}

\newcommand{\diag}{\mathsf{diag}\,}

\newcommand{\trace}{\mathsf{trace}\,}

\include{macro_list}

\usepackage[english]{babel}
\usepackage{graphicx}
\usepackage{epstopdf}

\usepackage[normalem]{ulem}

\newcommand{\krylov}[1]{\mathcal{K}_{#1}}

\newcommand{\TheTitle}{Efficient generalized Golub-Kahan based methods for dynamic inverse problems}

\title{{\TheTitle}}

\author{
  Julianne Chung
	\thanks{Department of Mathematics and Computational Modeling and Data Analytics Division, Academy of Integrated Science, Virginia Tech, Blacksburg, VA, \texttt{jmchung@vt.edu}}
  \and
  Arvind K. Saibaba
	\thanks{Department of Mathematics, North Carolina State University, Raleigh, NC, \texttt{asaibab@ncsu.edu}}
  \and
  Matthew Brown
	\thanks{Department of Mathematics, Virginia Tech, Blacksburg, VA, \texttt{brownm12@vt.edu} }
  \and
  Erik Westman
    \thanks{Department of Mining and Materials Engineering,  Virginia Tech, Blacksburg, VA, \texttt{ewestman@vt.edu} }
}

\newcommand{\post}{\boldsymbol{\Gamma}_{\text{post}}}
\newcommand{\posthat}{\widehat{\boldsymbol{\Gamma}}_{\text{post}}}

\begin{document}
\maketitle

\begin{abstract}
We consider efficient methods for computing solutions to and estimating uncertainties in dynamic inverse problems, where the parameters of interest may change during the measurement procedure. Compared to static inverse problems, incorporating prior information in both space and time in a Bayesian framework can become computationally intensive, in part, due to the large number of unknown parameters. In these problems, explicit computation of the square root and/or inverse of the prior covariance matrix is not possible.  In this work, we develop efficient, iterative, matrix-free methods based on the generalized Golub-Kahan bidiagonalization that allow automatic regularization parameter and variance estimation.   We demonstrate that these methods can be more flexible than standard methods and develop efficient implementations that can exploit structure in the prior, as well as possible structure in the forward model.  Numerical examples from photoacoustic tomography, deblurring, and passive seismic tomography demonstrate the range of applicability and effectiveness of the described approaches. Specifically, in passive seismic tomography, we demonstrate our approach on both synthetic and real data. To demonstrate the scalability of our algorithm, we solve a dynamic inverse problem with approximately $43,000$ measurements and $7.8$ million unknowns in under $40$ seconds on a standard desktop.
\end{abstract}

\vspace{2pc}

%


\section{Introduction}
The goal of an inverse problem is to use data, that is collected or measured, to estimate unknown parameters given some assumptions about the forward model~\cite{Vogel2002,Hansen2010}. In many applications, the problem is assumed to be static, in the sense that the underlying parameters do not change during the measurement process.  However, in many realistic scenarios such as in passive seismic tomography \cite{westman2012passive,zhang2009passive} or dynamic electrical impedance tomography \cite{schmitt2002efficient1,schmitt2002efficient2}, the underlying parameters of interest change dynamically. Incorporating prior information regarding temporal smoothness in reconstruction algorithms can lead to better reconstructions.  However, this presents a significant computational challenge since many large spatial reconstructions may need to be computed, e.g., at each time step and for many time points.  For example, in passive seismic tomography, geophones are used to collect measurements from seismic events (e.g., earthquakes) occurring in $1-2$ week intervals over $3-4$ months, and the goal is to obtain 3-dimensional spatial reconstructions of the elastic properties of the sub-surface for each time interval (e.g., to monitor changing stress conditions). As an other application, in medical imaging, during the data acquisition process the reconstruction algorithms need to account for patient motion and this can be modeled as a dynamic inverse problem.

Here we consider a discrete dynamic inverse problem with unknowns in space and time where
the goal is to reconstruct parameters $\bfs_i \in \mb{R}^{n_s}$ from observations $\bfd_i \in \mb{R}^{m_i}$ for $i = 1, \ldots, n_t$. Here $n_s$ refers to the number of spatial grid points, $n_t$ the number of time points, and $m=\sum_i m_i$ is the total number of measurements over all time points.  For some problems, the number of measurements may correspond to the number of sensors or spatial measurement locations and thus may be the same for all time points.  Let $$\bfs = \begin{bmatrix} \bfs_1 \\ \vdots\\ \bfs_{n_t}\end{bmatrix}\,,\qquad \mbox{and}\qquad
\bfd = \begin{bmatrix} \bfd_1 \\ \vdots\\ \bfd_{n_t}\end{bmatrix}\,,$$
then we are interested in the following dynamic inverse problem,
\begin{equation}
  \label{eq:invproblem}
 \bfd = \bfA \bfs + \bfvarepsilon \,,
\end{equation}
where $\bfA\in \bbR^{m \times n_s n_t}$ models the forward process which is assumed linear and $\bfvarepsilon$ represents noise or measurement errors in the data.  We assume that $\bfvarepsilon \sim \mc{N}(\bfzero,\bfR),$
where $\bfR$ is a positive definite matrix whose inverse and square root are inexpensive (e.g., a diagonal matrix with positive diagonal entries).
Given $\bfA$ and $\bfd,$ the goal of the inverse problem is to reconstruct $\bfs.$
Since these problems are typically ill-posed, regularization is often required to compute a reasonable solution.




To solve the dynamic inverse problem we adopt the Bayesian approach. In this approach, the measured data and the parameters to be recovered (here, the spacetime unknowns) are modeled as random variables. Additionally we assume that the prior distribution for $\bfs$ is modeled as a Gaussian distribution with mean $\bfmu = \begin{bmatrix} \bfmu_1\t & \cdots &  \bfmu_{n_t}\t\end{bmatrix}\t$ and positive-definite covariance matrix $\bfQ.$  That is, $\bfs \in \mc{N}(\bfmu, \lambda^{-2}\bfQ)$, where $\lambda$ is a (yet to be determined) scaling parameter for the precision matrix.  For dynamic problems, covariance matrices $\bfR$ and $\bfQ$ contain information in both space and time. Then Bayes' rule is used to combine the likelihood and the prior distribution and the posterior distribution,
\begin{align} \nonumber
 \pi(\bfs \mid\bfd) \propto & \> \pi(\bfd \mid\bfs) \pi(\bfs) \\ \label{e_bayes}
 \propto &\> \exp\left( -\frac12 \norm{\bfA \bfs -\bfd}{\bfR^{-1}}^2 - \frac{\lambda^2}{2} \norm{\bfs-\bfmu }{\bfQ^{-1}}^2\right),
\end{align}
where $\norm{\bfx}{\bfM} = \sqrt{\bfx\t \bfM \bfx} $ for any symmetric positive definite matrix $\bfM$.
The maximum a posteriori (MAP) estimate, which is the peak of the posterior distribution, can be obtained by minimizing the negative log likelihood of~\eqref{e_bayes}, i.e.,
\begin{equation}
  \label{eq:MAP}
 \bfs(\lambda) = \argmin_{\bfs \in \mb{R}^{n_s n_t}} \> \frac{1}{2}\|\bfA \bfs -\bfd \|_{\bfR^{-1}}^2 +  \frac{\lambda^2}{2}\| \bfs-\bfmu\|_{\bfQ^{-1}}^2\,,
\end{equation}

Notice that for dynamic inverse problems, computing the MAP estimate requires solving for $n_s n_t$ unknowns.  The main challenge here is that for the applications under consideration, $n_s$ is typically $O(10^5-10^6)$ and $n_t$ is typically $O(10^2-10^3)$. The resulting prior covariance matrices have $(n_sn_t)^2$ or $O(10^{14}-10^{18})$ entries. Storing such covariance matrices is completely infeasible, much less performing computations with them. Clearly there is a need for developing specialized numerical methods for tackling the immense computational challenges arising from dynamic inverse problems. Our strategy is to use a combination of highly structured representations of prior information along with efficient numerical methods that can exploit these representations.

\paragraph{Overview of main contributions.}  In this work, we adopt a Bayesian framework for solving dynamic inverse problems.  We derive two efficient methods for computing MAP estimates, where the distinguishing features of our approach compared to previous methods are that we can incorporate a wide class of spatiotemporal priors, include time-dependent observation operators, and enable automatic regularization parameter selection. The resulting solvers are highly efficient and scalable to large problem sizes. In addition to the MAP estimate, we develop an efficient representation of the posterior covariance matrix using the generalized Golub-Kahan (gen-GK) bidiagonalization. This low-rank approximation can be used for uncertainty quantification, by estimating the variance of the distribution. Using several real-world imaging applications (including both simulated and real data), we show that our methods are well suited for a wide class of dynamic inverse problems and we demonstrate scalability of our algorithms.

\paragraph{Related work.} The literature on dynamic inverse problems is large, and it is not our intention to provide a detailed overview. We mention a few related approaches that are relevant to our work.

A popular approach for spacetime reconstructions is the use of Kalman filters and smoothers. However, textbook implementations of these methods can be prohibitively expensive. This is because they require the storage and computation of covariance matrices that scale as $\mc{O}(n_s^2)$. One approach {is to use} an efficient representation of the state covariance matrix, as a low-rank perturbation of an appropriately chosen matrix~\cite{paninski2010fast,pnevmatikakis2013fast}. Efficient computational techniques for the Kalman filter especially tailored to  the random-walk forecast model were proposed by~\cite{li2014kalman,saibaba2015fastk}. In Section~\ref{sub:examples}, we show how this forecast model fits within our framework.

Significant simplifications can be made if we assume that measurement errors are independent in time and reconstruct the parameters of interest only using the data available from the current time step.
However, several other authors, see for example Schmitt and collaborators \cite{schmitt2002efficient1,schmitt2002efficient2}, have emphasized the importance of including temporal priors in many practical applications.
They considered a total-variation type temporal smoothness prior with a simple spatial prior (the identity matrix) and showed that their approach achieved superior results in faster computational time than other statistical approaches such as Kalman smoothers.

Our approach is more general in that we allow for a variety of spatial priors where the resulting covariance matrices are dense, unwieldy, and only available via matrix-vector multiplication.  Furthermore, we consider more general forward models and consider hybrid iterative approaches so that the regularization parameter $\lambda$ can be automatically estimated.  In previous studies, the regularization (or precision) parameter decoupled in space and time, and the resulting two parameters were required algorithmic inputs \cite{schmitt2002efficient1}.

In Section~\ref{sec:background}, we describe the problem set-up and address various changes of variables that can be used.  We also provide a brief overview of generalized hybrid methods to efficiently solve static inverse problems.  Efficient methods for approximating the MAP estimate, i.e., solving~\eqref{eq:MAP}, in the space-time formulation will be described in Section~\ref{sec:spacetime}, where special cases of problem structure will be considered for efficiency.  Efficient variance estimation methods based on the gen-GK bidiagonalization are described in Section~\ref{sec:sampling}.  Numerical results are presented in Section~\ref{sec:numerics} for various simulated imaging problems and for real data from passive seismic tomography.  Conclusions and discussions are provided in Section~\ref{sec:conclusions}.

\section{Problem set up and background}
\label{sec:background}

One main goal in the Bayesian framework is to efficiently compute the MAP estimate, and in this paper, we are mainly interested in cases where $\bfQ$ is a very large, dense matrix, so computing $\bfQ$ and $\bfQ^{-1}$ or their factorizations is not feasible.  Such scenarios arise, for example, when working with Gaussian random fields, in which case forming $\bfQ$ explicitly may not be feasible, but computing matrix-vector products (MVPs) with $\bfQ$ can be done efficiently \cite{saibaba2015fastc}.

First, we describe various problem formulations for computing the MAP estimate, and describe a change of variables so that gen-GK methods can be used.
Then, a brief overview of generalized hybrid methods is provided in Section~\ref{ss_genhybr} for completeness, and a discussion on various choices for modeling temporal priors is provided in Section~\ref{sub:modelingtemp}.   Connections to previous works are addressed in  Section~\ref{sub:examples}.


\subsection{Problem formulations}
Notice that the desired MAP estimate $\bfs(\lambda)$ is the solution to the system of equations
\begin{equation}\label{e_map}
(\bfA^\top\bfR^{-1}\bfA + \lambda^2 \bfQ^{-1})\bfs =\bfA^\top\bfR^{-1}\bfd + \lambda^2 \bfQ^{-1}\bfmu\,.
\end{equation}
For problems where matrix factorizations of $\bfR$ and $\bfQ$ are possible, a common approach to compute the MAP estimate is to solve the following \emph{general-form Tikhonov problem},
\begin{equation}
	\label{eq:genTik}
	\min_{\bfs} \> \frac{1}{2}\|\bfL_\bfR(\bfA \bfs -\bfd)\|_2^2 +  \frac{\lambda^2}{2}\| \bfL_\bfQ(\bfs-\bfmu)\|_2^2\,,
	\end{equation}
where $\bfQ^{-1} = \bfL_\bfQ\t \bfL_\bfQ$ and $\bfR^{-1} = \bfL_\bfR\t\bfL_\bfR$.
By transforming the problem to standard form, we get the \emph{priorconditioned problem} \cite{calvetti2015bayes,calvetti2005priorconditioners},
\begin{equation}
	\label{eqn:priorcondition}
	\min_{\bfx \in \mb{R}^{n_sn_t}} \> \frac{1}{2} \normtwo{\bfL_\bfR(\bfA\bfL_\bfQ^{-1}\bfx-\bfb)}^2 + \frac{\lambda^2}{2}\normtwo{\bfx}^2\,,
\end{equation}
where $\bfs = \bfmu+ \bfL_\bfQ^{-1} \bfx$.


In the applications that we consider, the covariance matrices can be very large and dense,  so the storage and computational costs to obtain factorizations {and/}or inverses {of $\bfQ$} can be prohibitive.
In order to avoid matrix factorizations of $\bfQ$ and/or expensive linear solves with $\bfQ$, a different change of variables was proposed in \cite{chungsaibaba2017}, where
\begin{equation}
  \label{eq:changevar}
\bfx \>\leftarrow  \> \bfQ^{-1}(\bfs-\bfmu), \qquad
 \bfb \> \leftarrow  \> \bfd - \bfA\bfmu\,,
\end{equation}
so that~\eqref{e_map} reduces to the modified system of equations
\begin{equation}\label{eqn:normal2}
(\bfA\t \bfR^{-1} \bfA \bfQ + \lambda^2 \bfI ) \bfx = \bfA\t \bfR^{-1}\bfb\,.
\end{equation}
In summary, with this change of variables, the MAP estimate is given by $\bfs(\lambda) = \bfmu + \bfQ \bfx(\lambda)$, where $\bfx(\lambda)$ is the solution to the following optimization problem
\begin{equation}
	\label{eqn:LS_tik_x}
	\min_{\bfx  \in \mb{R}^{n_sn_t}} \> \frac{1}{2}\|\bfA \bfQ \bfx -\bfb\|^2_{\bfR^{-1}} +\frac{ \lambda^2}{2}\|\bfx\|^2_\bfQ\,.
	\end{equation}
Hybrid iterative methods for  approximating $\bfx(\lambda)$ were described in \cite{chungsaibaba2017} for generic $\bfA, \bfR$, and $\bfQ$.  For completeness, we give a brief description of these methods in Section~\ref{ss_genhybr} and refer the interested reader to the paper for more details.

In this paper, we choose to focus on iterative methods; however, we briefly mention an alternative formulation. For relatively small problems where the number of overall measurements $m$ is small, i.e., $O(10^3-10^4)$, a simple change of variables along with the Sherman-Morrison formula can be used to avoid $\bfQ^{-1}$.  The MAP estimate can be computed as $\bfs(\lambda) = \bfmu + \bfQ\bfA\t\bfxi(\lambda)$ where 
\begin{equation}\label{e_xi}
 (\bfA \bfQ \bfA^\top + \lambda^2\bfR )\bfxi (\lambda) =  \bfd-\bfA\bfmu \,.
\end{equation}
Notice that in terms of solving linear systems, the number of unknowns has reduced from $n_sn_t$ to $m$.  A direct solver could be used to solve~\eqref{e_xi}, but forming $\bfQ\bfA^\top$ may be computationally prohibitive, costing $\mc{O}(n_s n_t m^2$).  Iterative methods could be used, but in these cases, it may be difficult to know a good regularization parameter a priori.  Further simplifications would be possible for problems where $m_i$ is constant for all $i=1,\dots,n_t$ and $\bfA$ is also a Kronecker product, i.e., $\bfA = \bfA_t \kron \bfA_s$ so that
\begin{align*}
 \bfQ \bfA^\top = & \>  (\bfQ_{t} \otimes \bfQ_s) (\bfA_t\t \otimes \bfA_s^\top) = \bfQ_t \bfA_t\t \otimes \bfQ_s \bfA_s^\top
\quad \mbox{and} \\
\bfA \bfQ \bfA^\top = & \>  (\bfA_t \bfQ_t \bfA_t\t) \otimes (\bfA_s\bfQ_s \bfA_s^\top)\,.
\end{align*}
As mentioned earlier, we do not pursue this approach since it is computationally expensive when the number of measurements is large.

\subsection{Generalized hybrid methods}\label{ss_genhybr}
Given matrices $\bfA$, $\bfR$, $\bfQ$, and vector $\bfb,$ with initializations $\beta_1 = \norm{\bfb}{\bfR^{-1}}, \bfu_1 = \bfb/\beta_1$ and $\alpha_1 \bfv_1 = \bfA\t \bfR^{-1} \bfu_1$, the $k$th iteration of the gen-GK bidiagonalization procedure generates vectors $\bfu_{k+1}$ and $\bfv_{k+1}$ such that
\begin{align*}
	\beta_{k+1} \bfu_{k+1} & = \bfA \bfQ \bfv_k -\alpha_k \bfu_k\\
		\alpha_{k+1} \bfv_{k+1} & = \bfA\t \bfR^{-1} \bfu_{k+1} -\beta_{k+1} \bfv_k,
	\end{align*}
where scalars $\alpha_i, \beta_i \geq 0$ are chosen such that $\norm{\bfu_i}{\bfR^{-1}} = \norm{\bfv_i}{\bfQ} = 1$. At the end of $k$ steps, we have
\[ \bfB_k \equiv \> \begin{bmatrix}
\alpha_1 \\ \beta_2 & \alpha_2 \\ & \beta_3 & \ddots \\ & & \ddots & \alpha_k \\ & & & \beta_{k+1}
\end{bmatrix}\,,  \qquad \bfU_{k+1} \equiv [\bfu_1,\dots,\bfu_{k+1}],\quad \mbox{and} \quad \bfV_k \equiv [\bfv_1,\dots,\bfv_k],\]
where the following relations hold up to machine precision,
\begin{align}\label{e_bk}
\bfU_{k+1}\beta_1 \bfe_1  =  &\> \bfb \\ \label{e_vk}
\bfA \bfQ \bfV_k = & \>\bfU_{k+1} \bfB_k \\ \label{e_uk}
\bfA\t \bfR^{-1} \bfU_{k+1} = & \> \bfV_k \bfB_k\t + \alpha_{k+1}\bfv_{k+1}\bfe_{k+1}\t\,.
\end{align}
Furthermore, in exact arithmetic, matrices $\bfU_{k+1}$ and $\bfV_k$ satisfy the following orthogonality conditions
\begin{equation}
	\label{eq:orthog}
	\bfU_{k+1}\t \bfR^{-1} \bfU_{k+1} = \bfI_{k+1} \qquad \mbox{and} \qquad \bfV_k\t \bfQ \bfV_k = \bfI_k.
\end{equation}

An algorithm for the gen-GK bidiagonalization process is provided in Algorithm~\ref{alg:wlsqr}. In addition to MVPs with $\bfA$ and $\bfA\t$ that are required for the standard GK bidiagonalization \cite{GoKa65}, each iteration of gen-GK bidiagonalization requires two MVPs with $\bfQ$ and two solves with $\bfR$ (which are assumed to be cheap); in particular, we emphasize that Algorithm~\ref{alg:wlsqr} avoids $\bfQ^{-1}$ and $\bfL_\bfQ$, due to the change of variables in~\eqref{eq:changevar}.

\begin{algorithm}[!h]
\begin{algorithmic}[1]
\REQUIRE Matrices $\bfA$, $\bfR$ and $\bfQ$, and vector $\bfb$.
\STATE $\beta_1 \bfu_1 = \bfb,$ where $\beta_1 = \norm{\bfb}{\bfR^{-1}}$
\STATE $\alpha_1 \bfv_1 = \bfA\t \bfR^{-1}\bfu_1$
\FOR {i=1, \dots, k}
\STATE $\beta_{i+1}\bfu_{i+1} = \bfA\bfQ\bfv_i - \alpha_i \bfu_i$, where $\beta_{i+1} = \norm{\bfA\bfQ\bfv_i - \alpha_i \bfu_i}{\bfR^{-1}}$
\STATE $\alpha_{i+1}\bfv_{i+1} = \bfA\t \bfR^{-1} \bfu_{i+1} - \beta_{i+1} \bfv_i$, where $\alpha_{i+1} = \norm{\bfA\t \bfR^{-1} \bfu_{i+1} - \beta_{i+1} \bfv_i}{\bfQ}$
\ENDFOR
\end{algorithmic}
\caption{generalized Golub-Kahan (gen-GK) bidiagonalization}
\label{alg:wlsqr}
\end{algorithm}

We seek solutions of the form $\bfx_{k}  = \bfV_k \bfz_{k} $, so that
\[ \bfx_k \in \text{Span}\{\bfV_k\} = \krylov{k}(\bfA\t\bfR^{-1}\bfA\bfQ,\bfA\t\bfR^{-1} \bfb) \equiv \mc{S}_k. \]
Define the residual at step $k$ as $\bfr_k \equiv \>  \bfA\bfQ\bfx_k-\bfb$. It follows from Equations~\eqref{e_bk}-\eqref{e_uk} that
{\[ \bfr_k \equiv \>  \bfA\bfQ\bfx_k-\bfb = \bfU_{k+1} \left( \bfB_k \bfz_k - \beta_1 \bfe_1\right)\,. \]}
To obtain coefficients $\bfz_k$, we take $\bfz_k$ that minimizes the genLSQR problem,
\begin{equation}\label{e_wlsqr}\min_{\bfx_k \in \mc{S}_k } \>\frac{1}{2}\norm{\bfr_{k}}{\bfR^{-1}}^2+\frac{\lambda^2}{2}\norm{\bfx_{k}}{\bfQ}^2\quad \Leftrightarrow \quad \min_{\bfz_k \in \mb{\bfR}^k} \> \frac{1}{2}\normtwo{ \bfB_k \bfz_k-\beta_1 \bfe_1 }^2 + \frac{\lambda^2}{2} \normtwo{\bfz_k}^2,
\end{equation}
where the gen-GK relations were used to obtain the equivalence.  Variants of this formulation, e.g., for LSMR, could be used as well \cite{chungsaibaba2017}. After computing a solution to the projected problem, an approximate MAP estimate can be recovered by undoing the change of variables,
\begin{equation}\label{eqn:undo_change}
\bfs_{k} =  \bfmu + \bfQ \bfx_{k}  = \bfmu + \bfQ\bfV_k\bfz_{k}\,,
\end{equation}
where, now, $\bfs_{k} \in \bfmu + \bfQ\mc{S}_k$.  For fixed $\lambda$ and in exact arithmetic, iterates of the genLSQR approach are mathematically equivalent to some pre-existing solvers (e.g., filtered GSVD solutions and priorconditioned solutions) \cite{chungsaibaba2017}.  However, if $\lambda$ is not known a priori, hybrid methods can take advantage of the shift-invariance property of Krylov subspaces and select $\lambda$ adaptively and automatically by utilizing well-known {SVD-based} regularization parameter selection schemes \cite{Hansen2010,HaHa93,Vogel2002} for the projected problem~\eqref{e_wlsqr}, {since} $\bfB_k$ is only of size $(k+1) \times k$. Henceforth, we refer to this approach as genHyBR.
Previous work on parameter selection within hybrid methods include \cite{ChNaOLe08,chung2015hybrid,Kilmer2001,gazzola2015krylov,renaut2010regularization,renaut2015hybrid}.

Although a wide range of regularization parameter selection methods can be used in our framework, in this paper we consider a variant of the generalized cross validation (GCV) approach. The GCV parameter is selected to minimize the GCV function \cite{GoHeWa79} corresponding to the general-form Tikhonov problem~\eqref{eq:genTik},
\begin{equation}\label{e_gcv}
G(\lambda) = \> \frac{n \norm{ \bfA \bfs(\lambda) - \bfd}{\bfR^{-1}}^2}{\left[\trace(\bfI_m-\bfL_\bfR \bfA \bfA_\lambda^{\dagger})\right]^2}\,,
\end{equation}
where $\bfA_\lambda^{\dagger} = (\bfA\t\bfR^{-1}\bfA + \lambda^2 \bfQ^{-1})^{-1}\bfA\t \bfL_\bfR\t$. We have assumed $\bfmu=\bfzero$ for simplicity.  At the $k$th iteration, the GCV parameter corresponding to the projected problem~\eqref{e_wlsqr} minimizes,
\begin{equation}\label{e_gcv_proj}
 G_\text{proj}(\lambda) \equiv \> \frac{k \normtwo{(\bfI - \bfB_k \bfB_{k,\lambda}^\dagger)\beta_1\bfe_1}^2}{\left[\trace(\bfI_{k+1} - \bfB_k\bfB_{k,\lambda}^\dagger)\right]^2 },
\end{equation}
where $\bfB_{k,\lambda}^\dagger = (\bfB_k\t\bfB_k+\lambda^2 \bfI)^{-1} \bfB_k\t$.
A weighted-GCV (WGCV) approach \cite{ChNaOLe08} has been suggested for use within hybrid methods, where a weighting parameter is introduced in the denominator of~\eqref{e_gcv_proj}.  We denote $\lambda_{\rm wgcv}$ to be the regularization parameter computed using WGCV.
As a benchmark for simulated experiments, we also consider the optimal regularization parameter $\lambda_{\rm opt}$, which minimizes the $2$-norm of the error between the reconstruction and the truth.

\subsection{Modeling prior covariances}
\label{sub:modelingtemp}
Following the geostatistical approach, we model the unknown field $s({\bfp},t)$ as a realization of a spatio-temporal random function $Z(\bfp,t)$ for $\bfp\in\bbR^d$ and $t \in \bbR$. We assume that the covariance function is stationary in space and stationary in time; for simplicity here, we also assume that the mean is zero. In other words, we assume that the covariance function satisfies \[ \text{cov}\{Z(\bfp_1,t_1),Z(\bfp_2,t_2)\} = C(\bfp_1-\bfp_2,t_1-t_2),\]
where  $C: \mb{R}^d \times \mb{R} \rightarrow \mb{R}_+$ is a positive definite covariance function.

We briefly review various formulations for spacetime covariance kernels, before delving into the specific choices of kernels we make. Perhaps the most convenient representation can be obtained if we make the assumption that the covariance function is separable in space and time, and isotropic in these variables, then $C$ takes the form
\[ C(\bfp,t) = C_S(\|\bfp\|)C_T(|t|) \qquad \forall (\bfp,t) \in \mb{R}^d \times \mb{R}, \]
where $C_S(\|\bfp\|)$ and $C_T(|t|)$ are isotropic, purely spatial and purely temporal
covariance functions, respectively. It can be readily seen that the resulting matrices $\bfQ$ have the Kronecker product structure. The separability assumption is common in the statistics literature~\cite{kyriakidis1999geostatistical,genton2007separable}, and tests for separability can be found in~\cite{fuentes2006testing}. The Kronecker product structure has computational advantages which we will exploit in Section~\ref{sec:spacetime} to develop efficient algorithms. While mathematically and computationally convenient, it is important to recognize the potential shortcomings of the separability assumption. The key issue is that the prior models do not allow for interactions in variability of space and time; for a detailed discussion, see~\cite{kyriakidis1999geostatistical}. 
However, even if the covariance kernel is not separable, one may approximate it using a separable covariance kernel~\cite{genton2007separable}, where the resulting covariance matrix approximation can be represented as a Kronecker product or a sum of Kronecker products.  Such approximations were studied in~\cite{loan1992approximation} and have been shown to be successful in the context of image deblurring, e.g.,~\cite{kamm1998kronecker,NaNgPe04,chung2015framework}.

Many nonseparable covariance kernels have been proposed that model space-time interactions of variability. One such approach uses $C(\bfp,t) = \varphi(\sqrt{c_1\|\bfp\|^2 + c_2|t|^2})$, where $c_1,c_2$ are weights that control the correlation of the space and time variables, and $\varphi(\cdot)$ is an appropriate covariance kernel. Another approach is to use a product-sum model
\[ C(\bfp,t) = a_0 C_S^0(\|\bfp\|)C_T^0(|t|) + a_1 C_S^1(\|\bfp\|) + a_2C_T^2(|t|)  \qquad \forall (\bfp,t) \in \mb{R}^d \times \mb{R}.\]
where $a_0$, $a_1$ and $a_2$ are nonnegative coefficients and $C_S^0$, $C_S^1$ and $C_T^0$, $C_T^2$ are isotropic, purely spatial and purely temporal covariance functions, respectively. A review of these covariance kernels is provided in~\cite{gneiting2006geostatistical}.

A wide variety of choices for $\bfQ$ can be included in our framework and thus incorporated in the methods described below. Next we give a few examples of temporal and spatial priors that are well-suited for our problems.

\paragraph{Specific choices of covariance kernels.} A common approach is to use Gaussian random fields where the entries of the covariance matrix are computed directly as
 $(\bfQ_t)_{ij} = \kappa(|t_i-t_j|)$, where $\{t_i\}_{i=1}^{n_t}$ are the time points. A popular choice for $\kappa(\cdot)$ is from the Mat\'{e}rn family of covariance kernels~\cite{rasmussen2006gaussian}, which form an isotropic, stationary, positive-definite class of covariance kernels. We define the covariance kernel in the Mat\'{e}rn class as
\begin{equation}\label{eqn:maternfamily}
   C_{\nu,\ell}(r) =   \frac{1}{2^{\nu-1}\Gamma(\nu)} \left(\frac{r\sqrt{2\nu}}{\ell}\right)^\nu K_\nu\left( \frac{r\sqrt{2\nu}}{\ell}\right)
\end{equation}
where  $\Gamma$ is the Gamma function, $K_\nu(\cdot)$ is the modified Bessel function of the second kind of order $\nu$, and $\ell$ is a scaling factor. The choice of parameter $\nu$ in equation~\eqref{eqn:maternfamily} defines a special form for the covariance. For example, when $\nu=1/2$, $C_{\nu,\ell}$ corresponds to the exponential covariance function, and if $\nu = 1/2+p$ where $p$ is a non-negative integer, $C_{\nu,\ell}$ is the product of an exponential covariance and a polynomial of order $p$. Also, in the limit as $\nu\rightarrow\infty$, $C_{\nu,\ell}$ converges to the Gaussian covariance kernel, for an appropriate scaling of $\ell$.
Another related family of covariance kernels is the $\gamma$-exponential function~\cite{rasmussen2006gaussian},
\begin{equation}
	\label{eq:gammaexp}
\kappa(r) = \exp\left(-(r/\ell)^\gamma\right) \qquad 0< \gamma \leq 2.
\end{equation}

Just as $C_T(\cdot)$ can be a $\gamma$-exponential function, or chosen from the Mat\'ern class, $C_S(\cdot)$ can also be chosen in the same way. Therefore, $\bfQ_s$ has entries $(\bfQ_s)_{ij} = \kappa(\|\bfp_i-\bfp_j\|)$, where $\{\bfp_i\}_{i=1}^{n_s}$ are the spatial locations. In the applications of interest, the number of spatial locations $n_s$ is much larger than the number of time points $n_t$; thus the storage of $\bfQ_s$ is challenging, as is employing it in iterative methods, since the cost of an MVP is $\mc{O}(n_s^2)$. Both the storage and computational cost can be reduced to $\mc{O}(n_s\log n_s)$ using the FFT based approach or $\mc{H}$-matrix approach. This has been reviewed in~\cite{saibaba2012efficient,saibaba2012application}.

\subsection{Other examples that fit our framework}
\label{sub:examples}
As mentioned in the introduction, the random-walk forecast model~\cite{vauhkonen1998kalman,kim2001image,soleimani2007dynamic,nenna2011application} was previously considered for its computational advantages~\cite{li2014kalman,saibaba2015fastk}. We show how this model also fits within our framework. {Assume that the state $\bfs_i$ undergoes the following dynamics for $i=1,\dots,n_t-1$}
 \[ \bfs_{i+1} = \bfs_i + \bfepsilon_i \qquad \bfepsilon_i \sim \mc{N}(\bf0,\bfQ_s),\]
with initial conditions $\bfs_1 \sim \mc{N}(\bf0,\bfQ_s)$. We can then express the distribution of the state $\bfs$ as
\[ \pi(\bfs) \propto \exp\left( -\frac12 \sum_{i=1}^{n_t-1} (\bfs_{i+1} -\bfs_i)\t\bfQ_s^{-1}(\bfs_{i+1} -\bfs_i) - \frac12 \bfs_1\t\bfQ_s^{-1}\bfs_1 \right).\]
{Thus, $\bfs$ is a Gaussian distribution with zero mean and  covariance matrix $\bfQ_t \kron \bfQ_s$, where}
\[ \bfQ_t^{-1} = \bmat{ 2 & -1 \\ -1 & 2 & -1 \\ & \ddots & \ddots & \ddots \\ & & -1 & 2 & -1 \\ &  & & -1& 1}. \]
The matrix $\bfQ_t$ has an explicit representation and is the so-called \verb|minij| matrix with $(i,j)$-th entry of $\bfQ_t$ equal to $\min\{i,j\}$. A similar representation is also available for the forecast model
${\bfs_{i+1} = \bar{\alpha_i} \bfs_i + \bar{\beta_i}\bfepsilon_i}$, but will not be considered here.

Another approach assumes that $\bfQ_t = (\bfL_t \t \bfL_t + \gamma \bfI)^{-1}$ where $\bfL_t$ is a sparse discretization of a differential operator and $\gamma$ is a small positive parameter to ensure $\bfQ_t$ is positive definite.  For example,
a common choice is to enforce smoothness in time by selecting
$$\bfL_t = \begin{bmatrix} \frac{1}{t_2-t_1} & - \frac{1}{t_2-t_1} & & &\\
 & \frac{1}{t_3-t_2} & - \frac{1}{t_3-t_2} & &\\
   &  & \ddots& \ddots &\\
   & & & \frac{1}{t_{n_t}-t_{n_t-1}} & - \frac{1}{t_{n_t}-t_{n_t-1}}
\end{bmatrix} \in \bbR^{(n_t-1) \times n_t}\,.$$
Although not derived within a Bayesian framework, Schmitt and collaborators \cite{schmitt2002efficient1,schmitt2002efficient2} considered such temporal priors along with a standard Tikhonov term to enforce spatial smoothness.  In fact, it is possible to show that their algorithm approximates the MAP estimate, which in our framework
corresponds to $\bfR= \bfI, \bfmu = \bfzero,$ and
\begin{equation}
\bfQ = \left(\lambda_s^2 \bfI + \lambda_t^2 \bfB\t \bfB\right)^{-1}
 = \underbrace{\left(\bfI + \frac{\lambda_t^2}{\lambda_s^2}\bfL_t\t \bfL_t\right)^{-1}}_{\bfQ_t} \kron \underbrace{\lambda_s^{-2}\bfI}_{\bfQ_s}
\end{equation}
where $\bfB = \bfL_t \kron \bfI.$
Here $\lambda_s$ and $\lambda_t$ correspond to  regularization parameters in space and time respectively, and the MAP estimate minimizes the function,
$$\Phi(\bfs) = \norm{\bfA\bfs - \bfd}{2}^2 + \lambda_s^2 \norm{\bfs}{2}^2 + \lambda_t^2 \norm{\bfB \bfs}{2}^2\,.$$

In this paper, we focus on Mat\'{e}rn kernels and their covariance matrices for both the spatial and temporal priors.



\section{Generalized hybrid methods for dynamic inverse problems}
\label{sec:spacetime}
In this section, we describe various approaches based on the gen-GK bidiagonalization for computing MAP estimates for dynamic inverse problems. These iterative approaches are desirable for problems where both $\bfA$ and $\bfQ$ may be extremely large, or for problems where these matrices are not explicitly stored but can be accessed via function calls to compute MVPs with $\bfA, \bfA\t,$ and $\bfQ$ efficiently. In Section~\ref{sub:simultaneous}, we describe an ``all-at-once'' generalized hybrid method that  requires MVPs with covariance matrix $\bfQ$.  This requirement is quite general and includes matrices such as $\bfQ$ being a Kronecker product, a sum of Kronecker products, or a convolution operator.
Then for problems where the number of time points is small and $\bfA$ and $\bfR$ are also Kronecker products, we describe an efficient decoupled approach in Section~\ref{sub:decoupled}.


\subsection{{S}imultaneous generalized hybrid approach}
\label{sub:simultaneous}
The first approach we consider for solving dynamic inverse problems is to use genHyBR as summarized in Section~\ref{ss_genhybr} to solve for all unknown variables (e.g., in space and time) simultaneously.
Since the number of unknowns can be quite large in the ``all-at-once'' approach and the gen-GK vectors $\bfV_k$ must be stored {for hybrid methods}, we assume that solutions can be captured in relatively few iterations or that appropriate preconditioning can be used so that $k$ remains small.
Next we describe efficiencies can can be gained for problems {where $\bfQ$ is a Kronecker product,}
 but we reiterate that the simultaneous approach has applicability beyond the cases presented here.

  %





For problems where $\bfQ$ is a Kronecker product~\cite{laub2005matrix}, MVPs with $\bfQ$ can be computed efficiently as
\[ \bfQ\bfx = (\bfQ_t\otimes \bfQ_s)\bfx = \mathsf{vec}(\bfQ_s\bfX\bfQ_t\t)\,,\]
where $\mathsf{vec}$ and $\mathsf{mat}$ are operations such that $\mathsf{vec}$ unfolds a matrix $\bfX \in \mb{R}^{m\times n}$ into a vector by stacking column-wise and $\mathsf{mat}$ folds it back, i.e.,
\[ \mathsf{vec}(\bfX)\in \mb{R}^{mn\times 1} \qquad \mathsf{mat}[\mathsf{vec}(\bfX)] = \bfX. \]
Assuming the cost of an MVP with $\bfQ_s$ is $\mc{O}(n_s\log n_s)$, the cost of $\bfQ\bfx$ is $\mc{O}(n_t n_s \log n_s + n_s n_t^2)$, which is significantly smaller than the naive cost of $\mc{O}(n_s^2n_t^2)$. Further reductions in computational cost can be achieved by parallelizing the matrix-matrix multiplications, e.g., MVPs of $\bfQ_s$ with the columns of $\bfX$.



\paragraph{Separable forward operator.}

If, additionally, the number of measurements at each timestep is the same, which we denote by $\bar{m}$, and $\bfA = \bfA_t \otimes \bfA_s$ and
 $\bfR = \bfR_t \otimes \bfR_s$, where $\bfR_{t}\in\bbR^{n_t \times n_t}$ and $\bfR_s \in \bbR^{\bar{m} \times \bar{m}}$ are the temporal and spatial noise covariance matrices respectively, then MVPs required for the gen-GK bidiagonalization algorithm can be computed efficiently.
 That is, for vectors $\bfx \in \mb{R}^{n_sn_t}$ and $\B{y} \in \mb{R}^{\bar{m}n_t}$
$$ \bfA \bfQ\bfx = \mathsf{vec}(\bfA_s\bfQ_s\bfX\bfQ_t\t \bfA_t^\top),
\quad \mbox{and} \quad \bfA^\top\bfR^{-1}\bfy =   \mathsf{vec}(\bfA_s^\top\bfR_s^{-1}\bfY \bfR_t^{-\top} \bfA_t),$$
where $\bfX \equiv \mathsf{mat}(\bfX) \in \mb{R}^{n_s\times n_t}$ and $\bfY \equiv \mathsf{mat}(\bfy)\in \mb{R}^{\bar{m} \times n_t}$.




\subsection{{D}ecoupled generalized hybrid approach}
\label{sub:decoupled}
For problems where $\bfA, \bfR$, and $\bfQ$ are all Kronecker products and $n_t$ is relatively small such that $\bfQ_t$ and its factor $\bfL_t$ are feasible, we describe a
decoupled genHyBR approach.
The normal equations corresponding to the weighted least-squares problem in~\eqref{eqn:LS_tik_x} can be written as
\begin{equation}
 (\bfQ\t \bfA\t \bfR^{-1} \bfA \bfQ + \lambda^2 \bfQ)\bfx = \bfQ\t \bfA\t \bfR^{-1} \bfb \,.
\end{equation}
We can exploit the fact that
$\bfQ = \bfL_t\t\bfL_t \otimes \bfQ_s \,,$
and let $\bfy = (\bfL_t \otimes \bfI)\bfx$, alternatively $\bfY = \bfX\bfL_t\t$. We introduce the variable $\bfH = \bfA\t \bfR^{-1} \bfA$; similarly, we also define $\bfH_s = \bfA\t_s \bfR^{-1}_s \bfA_s$ and $\bfH_t = \bfA_t\t \bfR_t^{-1} \bfA_t$. The normal equations simplify, when we left-multiply by $\bfL_t^{-\top} \otimes \bfI$ to obtain
\begin{align}
 [ (\bfL_t \otimes \bfQ_s) \bfH  (\bfL_t\t  \otimes \bfQ_s) + \lambda^2 (\bfI \otimes \bfQ_s)] \bfy & =  (\bfL_t \otimes \bfQ_s) \bfA\t \bfR^{-1} \bfb.
\end{align}
Using the properties of Kronecker products, this equation can alternatively be written as the following {g}eneralized-Sylvester equation
\[
\bfQ_s \bfH_s \bfQ_s \bfY \bfL_t \bfH_t \bfL_t\t + \lambda^2 \bfQ_s \bfY  = \bfC\widehat{\bfA}_t, \]
where, for simplicity, we introduce $\bfC \equiv \bfQ_s \bfA_s\t \bfR_s^{-1}\bfB\bfR_t^{-1/2} $ and $\widehat{\bfA}_t\equiv \bfR_t^{-1/2}\bfA_t \bfL_t\t$.

Let $\widehat{\bfA}_t  = \bfU_t \bfSigma_t \bfV_t\t$ be its singular value decomposition, then $\bfL_t \bfH_t \bfL_t\t = \bfV_t\bfSigma_t^2 \bfV\t_t$. We make another change of variables $\bfZ \equiv \bfY \bfV_t = \bfX \bfL_t\t \bfV_t,$ and we get
\begin{align}
\bfQ_s \bfH_s \bfQ_s \bfZ \bfSigma^2 + \lambda^2 \bfQ_s \bfZ & = \bfC \bfU_t \bfSigma_t,
\end{align}
where we have multiplied on the right by $\bfV_t$. Expanding the above expression column-wise, the key observation is that all the equations decouple so that
\begin{equation}(\sigma_i^2  \bfQ_s \bfH_s \bfQ_s + \lambda^2 \bfQ_s)  \bfz_i = \sigma_i \bfC  \bfu_i \quad i = 1, \cdots, n_t\,.
\end{equation}
Notice that for $\sigma_i = 0$, $\bfz_i = \bf0$, whereas for $\sigma_i>0,$ the solution can be obtained by solving the least-squares problem,
\begin{equation}\label{e_decoupled}
\min_{\bfz_i \in \mb{R}^{n_s}} \frac12 \| \sigma_i  \bfA_s \bfQ_s \bfz_i  - \bfB\bfR_t^{-1/2}\bfu_i \|_{\bfR_s^{-1}}^2 + \frac{\lambda^2}{2} \|\bfz_i\|_{\bfQ_s}^2\,,
\end{equation}
{which can be done using the}  genHyBR method.  Note that a transformation back of variables must be made. This is summarized in Algorithm~\ref{alg:decoupled}.

\begin{algorithm}[!ht]
\begin{algorithmic}[1]
\REQUIRE $\bfA_s, \bfA_t, \bfR_s, \bfR_t, \bfL_t, \bfB$
\STATE Compute the SVD of $\widehat{\bfA}_t = \bfR_t^{-1/2} \bfA_t \bfL_t\t  = \bfU_t \bfSigma_t \bfV_t\t$
\FOR{i=1, $\cdots$, $n_t$}
\STATE Apply Algorithm~\ref{alg:wlsqr} to~\eqref{e_decoupled} to obtain $\bfz_i$.
\ENDFOR
\STATE Form $\bfZ = \bmat{\bfz_1 & \dots & \bfz_{n_t}}$ and $\bfX = \bfZ\bfV_t\t\bfL_t^{-\top}$ and $\bfS = \bfQ_s \bfX \bfQ_t\t$
\end{algorithmic}
\caption{Decoupled genHyBR}
\label{alg:decoupled}

\end{algorithm}

%
%
This approach is embarrassingly parallel, which may lead to enormous computational savings.
 {Additionally, solving the sequence of decoupled systems may be done efficiently either by recycling of Krylov subspaces or by effective preconditioning.}

\section{Estimation of posterior variances using gen-GK}
\label{sec:sampling}
For the problem considered here, the posterior distribution $\bfs | \bfd$ is Gaussian with
\[ \bfs | \bfd \sim \mc{N}(\post \bfA\t\bfR^{-1}\bfb, \post) \qquad \post \equiv ( \lambda^2 \bfQ^{-1}+\bfH)^{-1},\]
where $\bfH = \bfA\t\bfR^{-1}\bfA$.  The MAP estimate {corresponds to} the mode of the posterior distribution {and} provides information regarding the ``most likely'' estimate. On the other hand, the posterior variance defined as the diagonals of the posterior covariance $\post$
provides a measure of the spread of the posterior distribution around the posterior mean.
For {static} inverse problems, estimating the posterior covariance matrix is known to be computationally challenging~\cite{saibaba2012efficient,saibaba2015fastc}.
For dynamic problems, the problem is further exacerbated since the posterior covariance matrix is of size $n_s n_t \times n_s n_t$. Moreover, this matrix is dense and forming it explicitly {to obtain the diagonal entries is} computationally infeasible. In this work, we use intermediate information from genHyBR for computing a MAP estimate to estimate the posterior variance. Following Section~\ref{sec:spacetime}, we consider a simultaneous and a decoupled approach. 

Before we explain how to compute variances, we make the following remark. In the sequel, we will use Algorithm~\ref{alg:wlsqr} with one minor modification, namely, at each step we explicitly re-orthogonalize the vectors in $\bfU_k$ and $\bfV_k$. Numerical experience suggests that this marginally increases the computational cost by $\mc{O}(k^2(m+n))$ but considerably improves the accuracy of the variance computations.

\subsection{Simultaneous approach}
Recall that after $k$ iterations of the gen-GK bidiagonalization process, we have matrices $\bfB_k,\bfU_k$ and $\bfV_k$ satisfying relations~\eqref{e_bk}-\eqref{e_uk} and~\eqref{eq:orthog}. Let $\bfB_{k}\t\bfB_{k} = \bfW_k\bfTheta_k\bfW_k\t$ be the eigenvalue decomposition with
{eigenvalues $\theta_1, \ldots, \theta_k$}
and let $\bfZ_k = \bfQ\bfV_k\bfW_k$, then we get the following low-rank approximation
\begin{equation}
  \label{eq:GKapprox}
  \bfQ\bfH\bfQ \approx \> \bfQ\bfV_k\bfB_{k}\t\bfB_{k} \bfV_k^\top\bfQ = \> \bfZ_k\bfTheta_k\bfZ_k\t\,.
\end{equation}

Using~\eqref{eq:GKapprox} and the Woodbury formula, we obtain the approximation
\begin{align*}
\post \approx & \>  \bfQ(\lambda^2 \bfQ + \bfZ_k \bfTheta_k\bfZ_k^T)^{-1}\bfQ \\
 = &  \>   {\lambda^{-2}} \bfQ - \lambda^{-2}\bfZ_k (\bfI_k  + \lambda^2\bfTheta_k^{-1})^{-1}\bfZ_k\t\\\
 = &  \>  {\lambda^{-2}} \bfQ - \bfZ_k\bfDelta_k\bfZ_k\t \equiv \posthat
\end{align*}
where
\[\bfDelta_k \equiv \lambda^{-2} \begin{bmatrix} \frac{\theta_1}{\theta_1 + \lambda^2} & & \\ & \ddots & \\ & & \frac{\theta_k}{\theta_k + \lambda^2} \end{bmatrix}  \in \mb{R}^{k\times k} .\]
Notice that we have an efficient representation of $\posthat$ as a low-rank perturbation of the prior $\lambda^{-2}\bfQ$.
In summary, diagonal entries of $\posthat$ can provide estimates of diagonal entries of $\post$, where the main computational requirement is to obtain the diagonals of $\bfQ$  
 and the diagonals of the rank-$k$ perturbation.  Therefore, the only additional computational cost for estimating the posterior variance is $\mc{O}(k^3 + k^2n_sn_t)$.

An approximation of this kind was previously explored in~\cite{saibaba2015fastc,flath2011fast,bui2012extreme,bui2013computational}; however, the error estimates developed in the above references assume that the exact eigenpairs are available. If the Ritz pairs converge to the exact eigenpairs of the matrix $\bfQ\bfH\bfQ$, then furthermore, the optimality result in~\cite[Theorem 2.3]{spantini2015optimal} applies here as well.

\subsection{Decoupled approach}
\label{sub:decvar}For cases where $\bfA=\bfA_t \kron \bfA_s,$ we develop a similar strategy for estimating the posterior variances by exploiting the decoupled structure described in Section~\ref{sub:decoupled}.
In contrast to the simultaneous approach, a different Krylov subspace is constructed for each time step.
Using notation defined in Section~\ref{sub:decoupled}, the posterior covariance matrix is given by
\begin{align*}
  \post = & \> ( \bfQ_t^{-1} \otimes \lambda^2 \bfQ_s^{-1} + \bfH_t \otimes \bfH_s )^{-1} \\
  = & \> (\bfL_t\t \bfV_t \otimes \bfI) (\bfI\otimes\lambda^2\bfQ_s^{-1} + \bfSigma_t^2 \otimes \bfH_s)^{-1}(\bfV_t\t\bfL_t \otimes \bfI)\,,
\end{align*}
where the matrix in the center is a block-diagonal matrix whose diagonal blocks are $(\lambda^2\bfQ_s^{-1} +\sigma_i^2 \bfH_s)^{-1}$ for $i = 1,\ldots, n_t\,.$  Analogous to the simultaneous approach, gen-GK approximations for each $i$ denoted by $\bfB_{k,i}$ and $\bfV_{k,i}$ can be used to get low-rank approximations
$$(\lambda^2\bfQ_s^{-1} +\sigma_i^2 \bfH_s)^{-1} \approx \lambda^{-2}\bfQ_s- \bfZ_{k,i}\bfDelta_{k,i}\bfZ_{k,i}\t , \qquad i=1,\dots,n_t $$
where $\bfB_{k,i}\t \bfB_{k,i} = \bfW_{k,i} \bfTheta_{k,i} \bfW_{k,i}\t$ is an eigenvalue decomposition,
$$\bfZ_{k,i} = \bfQ_s \bfV_{k,i} \bfW_{k,i},  \qquad \text{and}\qquad \bfDelta_{k,i} = \lambda^{-2}(\bfI_k+\lambda^2 \bfTheta_{k,i}^{-1})^{-1}.$$
Notice that depending on the stopping criteria for the gen-GK process, each time point may have a different rank $k$.   We omit this dependence for clarity of presentation, 
and denote the low-rank blocks $\bfD_i = \bfZ_{k,i}\bfDelta_{k,i}\bfZ_{k,i}\t$.

In summary, an approximation to the posterior covariance matrix is given by
\[ \post \approx  \lambda^{-2}\bfQ - (\bfL_t\t\bfV_t \otimes \bfI) \bmat{\bfD_1 \\ & \ddots \\ & &  \bfD_{n_t}} (\bfV_t\t\bfL_t \otimes \bfI),\]
where the diagonals can be computed for timestep $i=1,\dots,n_t$ as
\[ \diag\left[(\bfe_i\t \otimes \bfI)\post(\bfe_i\otimes \bfI)\right] \approx \lambda^{-2}([\bfQ_t]_{ii} \otimes \diag(\bfQ_s))  - \sum_{j=1}^{n_t} (\bfe_j\t \bfV_t\t \bfL_t \bfe_i)^2 \,\diag(\bfD_j) \,,\]
where $[\bfQ_t]_{ii}$ is the $i$th diagonal entry of $\bfQ_t$ and the operation $\diag(\cdot)$ returns a vector containing the diagonals of a matrix.

\section{Numerical Results}
\label{sec:numerics}
In this section, we provide three examples from image reconstruction.  The first is a model problem from dynamic photoacoustic tomography (PAT) reconstruction under motion that illustrates the significant impact of including a temporal prior.  Then we consider an image deblurring problem that can exploit the decoupled framework of Section~\ref{sub:decoupled}, and we finish with a severely ill-posed passive seismic tomography (PST) problem, where we include results from real field measurements.  In all of the results, genHyBR solutions correspond to~\eqref{eqn:undo_change} where $\bfx_k$ is the solution to~\eqref{e_wlsqr}.

\subsection{Dynamic photoacoustic tomography (PAT) under motion}
PAT is a hybrid imaging modality that combines the rich contrast of optical imaging with the high resolution of ultrasound imaging, thereby producing higher resolution in-vivo images with lower patient risk (e.g., requiring no ionizing radiation and no contrast agents) and lower cost and inconvenience than other imaging modalities.  In modern PAT systems, transducers are rotated around an object and data is acquired in time. Most current methods for PAT image reconstruction such as explicit inversion formulas, time reversal, or series solution are not suited for such systems. A significant limitation in extending these methods to dynamic PAT is that accurate tomographic reconstruction relies on expensive motion estimation and artifact removal software, since the object being imaged may move during data acquisition. Chung and Nguyen \cite{chung2017motion} recently studied a continuous model of PAT reconstruction under motion for specific parameterized motion models and developed specialized algorithms for these scenarios. However, in this work we consider PAT reconstruction for more general motion deformations and seek a sequence of reconstructions, rather than just one reconstruction, by incorporating a temporal prior and more informative image priors.  We follow the literature in dynamic tomography (e.g., \cite{katsevich2011local,hahn2014efficient,hahn2015dynamic}), and note that since we do not consider parameterized models, our framework can incorporate other realistic scenarios (other than motion) such as non-stationary optical illumination and inaccurate transducer responses \cite{sheng2015constrained,zhang2008simultaneous,wang2011photoacoustic,lou2016impact}.

We consider the discrete problem where $\bfs_i\in \bbR^{n_s}$ is the discretized desired solution\footnote{Here we assume the 2D image is vectorized column-wise.} at time point $i$ and let $\left\{\bfz_i\right\}, i = 1,\ldots, n_t$ denote the locations of the transducers.   At each transducer location $\bfz_i$, assume there are $r$ radii and let $\bfA_i \in \bbR^{r \times N}$ be the corresponding projection matrix for that location (i.e., $\bfA_i \bfs_i$ is the discrete circular Radon transform of $\bfs$ on circles centered at $z_i$).
Thus, the observed spherical projection measurements for all $r$ radii are contained in vector
$$\bfd_i = \bfA_i \,\bfs_i + \bfe_i \,\, \in \bbR^{r}$$
where $\bfe_i \in \bbR^{r}$ is additive Gaussian noise that is independent and identically distributed.
Then the forward model has the form
\begin{equation}
  \label{eqn:PATforward}
  \begin{bmatrix} \bfd_1 \\ \vdots \\ \bfd_{n_t}
  \end{bmatrix} = \begin{bmatrix}
    \bfA_1 & & \\ &\ddots & \\&& \bfA_{n_t} \end{bmatrix} \begin{bmatrix} \bfs_1 \\ \vdots \\ \bfs_{n_t}\end{bmatrix} + \begin{bmatrix}\bfe_1 \\ \vdots \\ \bfe_{n_t}
  \end{bmatrix}
\end{equation}
where the goal is to estimate the desired images $\bfs_i$, given the observations ${\bfd_i.}$

\begin{figure}[ht]
 \centering
   \begin{tabular}{ccc}
     \includegraphics[width=.4\textwidth]{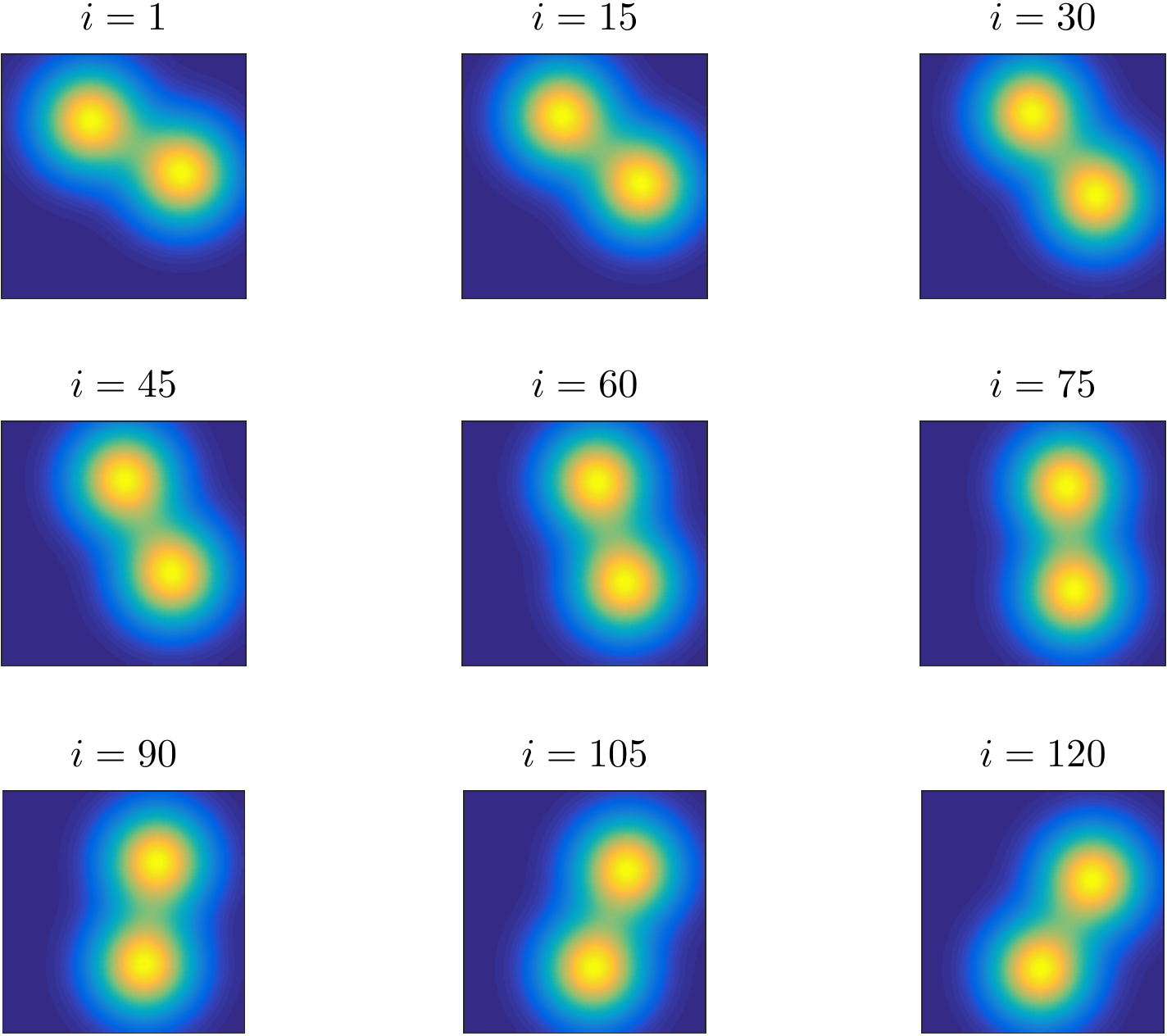} &
    \includegraphics[width=.15\textwidth]{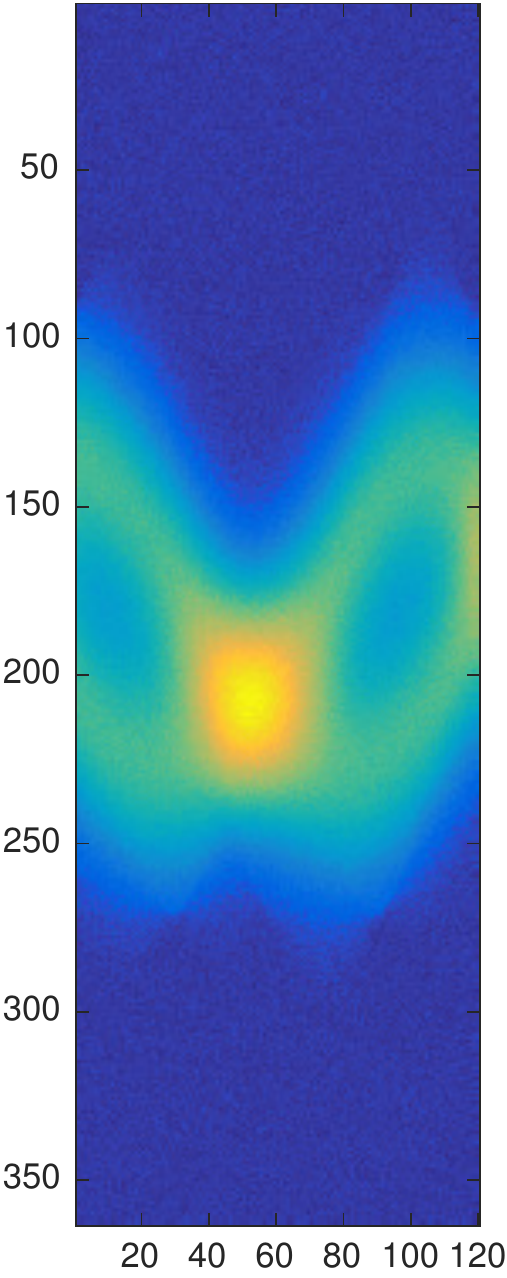} &
    \includegraphics[width=.25\textwidth]{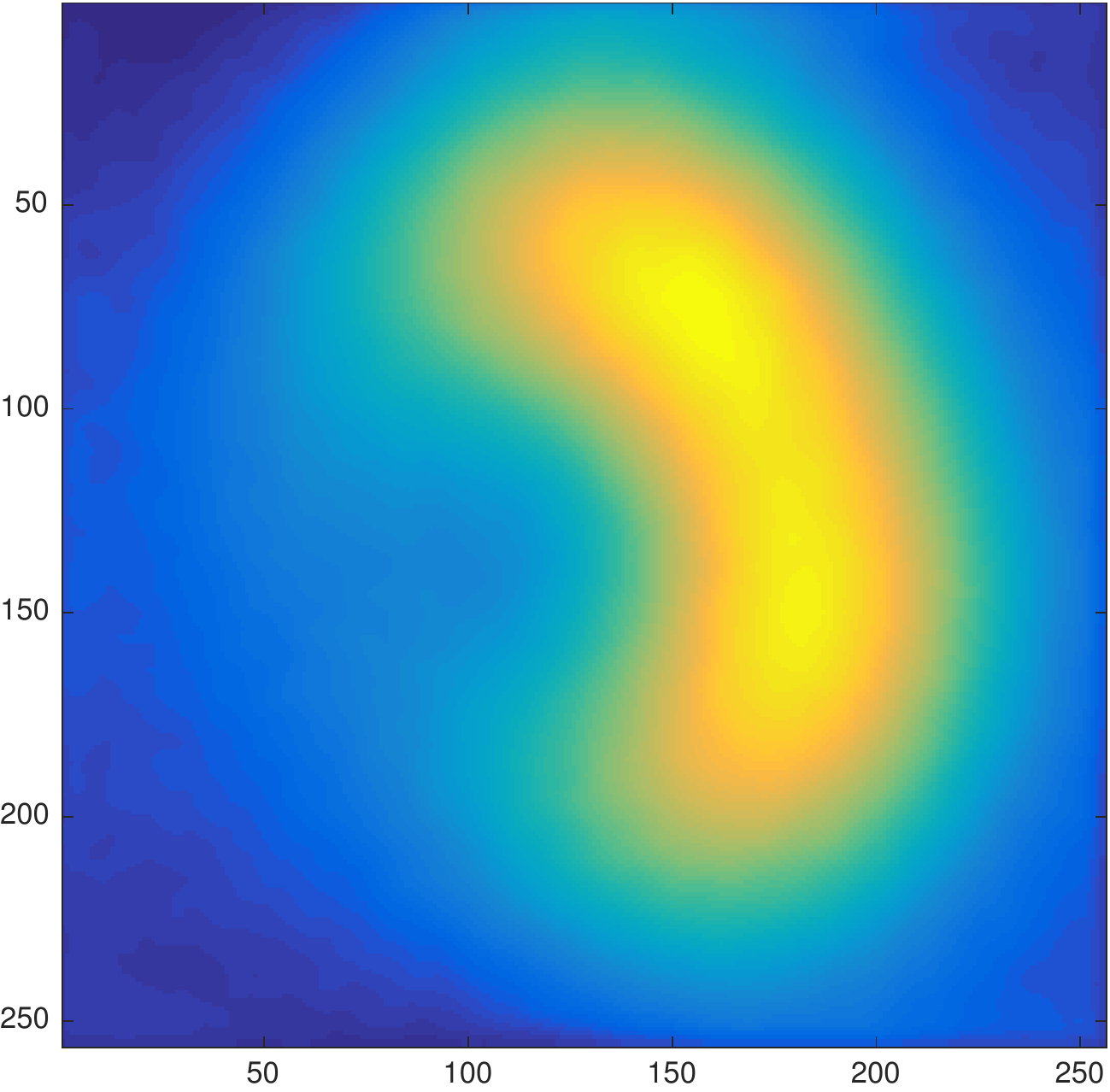}\\
  (a) Sample true images & (b) Sinogram & (c) Static reconstruction
   \end{tabular}
  \caption{PAT reconstruction problem.  Sample true images are provided in (a), observations in the form of a sinogram are provided in (b), and a static reconstruction (i.e., ignoring temporal changes) of the data is provided in (c) for comparison.}\label{fig:PAT_problem}
\end{figure}

For this example, we consider $120$ true images of size $256 \times 256$ that were generated with two Gaussians with fixed width, and rotating counterclockwise.  See Figure~\ref{fig:PAT_problem}(a) for sample true images. Measurements were taken at $120$ equidistant angles between $0$ and $357$ at $3$ degree intervals, and each projection consists of $363$ radii. White noise was added to the observations with $\bfR = .0082^2 \bfI$ (this corresponds to a noise level of $0.04$).  Results presented here use $\bfR$ as above, although noise estimation algorithms could be used \cite{Do95}.
The sinogram of size $363 \times 120$ is shown in Figure~\ref{fig:PAT_problem}(b), providing a total number of $43,560$ observations.

For the dynamic inverse problem, the number of unknowns is $256*256*120=7,864,320$.  
We provide, for comparison, a static reconstruction in Figure~\ref{fig:PAT_problem}(c), where we solve the following (inaccurate) model problem,
$$\min_{\bfs \in \mb{R}^{n_s}}\, \left\| \bfd -  \bfA\bfs \right\|_{\bfR^{-1}}^2 + \lambda^2 \norm{\bfs}{\bfQ_s^{-1}}^2 $$
with $\bfd = \begin{bmatrix} \bfd_1^\top& \cdots & \bfd_{n_t}^\top
\end{bmatrix}^\top$ and $\bfA = \begin{bmatrix}
  \bfA_1^\top & \cdots & \bfA_{n_t}^\top \end{bmatrix}^\top .$
The least squares problem above is solved using genHyBR where the regularization parameter is picked using the WGCV criterion to obtain one reconstructed image of size $256 \times 256$.  Here we set $\bfQ_s$ to be a covariance matrix that is determined from the Mat\'ern covariance function $C_S=C_{1,.01}(\cdot)$ as described in Section~\ref{sub:modelingtemp}.  It is evident that the static reconstruction is able to locate the object but can neither distinguish the objects nor provide dynamic information.

Next we consider three cases that use simultaneous genHyBR to solve~\eqref{eqn:PATforward}:
\begin{itemize}
  \item genHyBR where $\bfQ$ is generated from a Mat\'ern kernel $C_{1,.01}(\sqrt{c_1 \|\bfp\|^2 + c_2 |t|^2} )$ where $c_1 = 1$ and $c_2 = 0.0025$.
  Here, $\bfQ$ can not be represented as a Kronecker product, but MVPs can still be done efficiently.
  \item genHyBR with $\bfQ = \bfQ_t \kron \bfQ_s$ where $\bfQ_t = \bfI$ and $\bfQ_s$ corresponds to $C_S(\cdot)=C_{1, .01}(\cdot)$.
  \item genHyBR with $\bfQ = \bfQ_t \kron \bfQ_s$ where $\bfQ_t$ and $\bfQ_s$ correspond to $C_T(\cdot)=C_{\infty, .01}(\cdot)$ and $C_S(\cdot)=C_{1, .01}(\cdot)$ respectively.
\end{itemize}

\begin{figure}[ht]
 \centering
   \includegraphics[width=\textwidth]{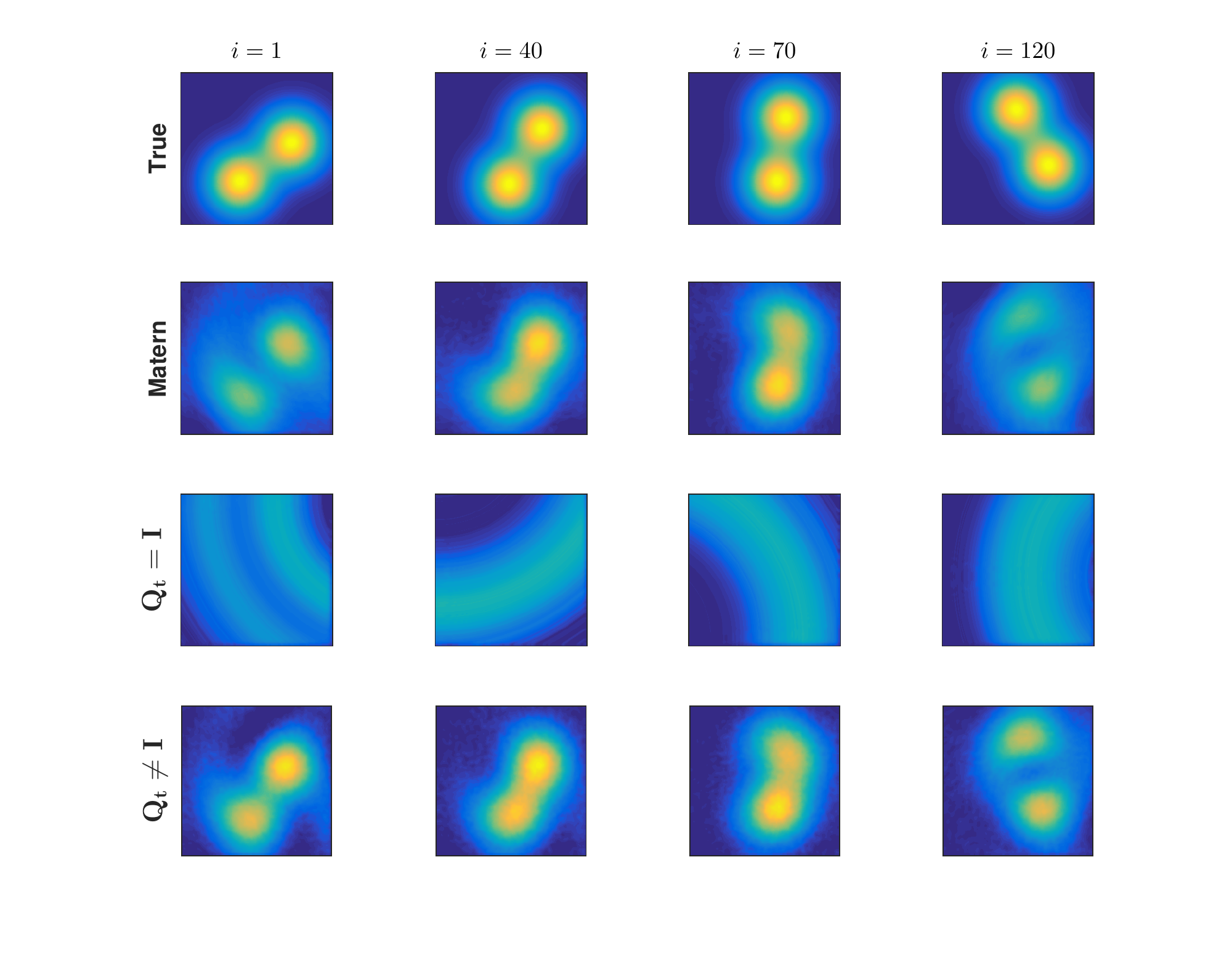}
  \caption{Dynamic PAT reconstructions using genHyBR with WGCV at various time points ($i=1, 40, 70, 120$).  These results show the importance of incorporating a temporal prior in dynamic PAT setups.}\label{fig:PAT}
\end{figure}
Image reconstructions for various time points, along with the corresponding true images, are provided in Figure~\ref{fig:PAT}, where all of the results use the WGCV parameter after $10$ iterations.  Compared to the static reconstruction in Figure~\ref{fig:PAT_problem}(c), the Mat\'ern reconstructions in the second row reveal changes over time.  However, the more striking comparison occurs when using the separable covariance functions and comparing the results with $\bfQ_t = \bfI$ to those with $\bfQ_t \neq \bfI$ (c.f.~rows 3 and 4 in Figure~\ref{fig:PAT}).  Dynamic PAT is a severely underdetermined problem, and this example illustrates that including a temporal prior can be crucial to revealing dynamics of the imaged object.  In Table~\ref{tab:PATregpar} we provide the computed regularization parameters for each approach, along with the relative errors computed as $\norm{\bfs_k-\bfs_{\true}}{2}/\norm{\bfs_\true}{2}$.  These values are consistent with the quality of the reconstructions in Figure~\ref{fig:PAT} {and are not significantly improved with reorthogonalization of gen-GK vectors}.  {For reconstructions that assume $\bfQ$ is a Kronecker product, genHyBR took around $37$ seconds and, with reorthogonalization, around $106$ seconds\footnote{{All timings were recorded on a MacPro, OSX Yosemite, 2.7 GHz 12-Core Intel Xeon E5, 64G memory in Matlab 2014b using default computing options.}}.  Partial reorthogonalization could be used but was not investigated here.}
\begin{table}[ht]
  \centering

  \begin{tabular}{|c|c|c|c|c|}\hline
    & $\lambda_{\rm wgcv}$ & relative error & without reorth (sec) & with reorth (sec)\\ \hline
    Mat\'ern & 40.32 & 2.9923e-01 & 84.5  & 273.7 \\ \hline
     $\bfQ_t = \bfI$ & 6.24 & 6.4575e-01 & 36.5  & 107.9 \\ \hline
     $\bfQ_t \neq \bfI$ & 35.93 & 2.3406e-01 & 36.7  & 105.4 \\ \hline
  \end{tabular}
\caption{Regularization parameters computed using WGCV and relative errors for PAT reconstructions {corresponding to reorthogonalization.  CPU time (in seconds) to obtain dynamic PAT reconstructions without and with reorthogonalization of gen-GK vectors.}}
\label{tab:PATregpar}
\end{table}

\begin{figure}[ht]
  \centering
    \includegraphics[width=\textwidth]{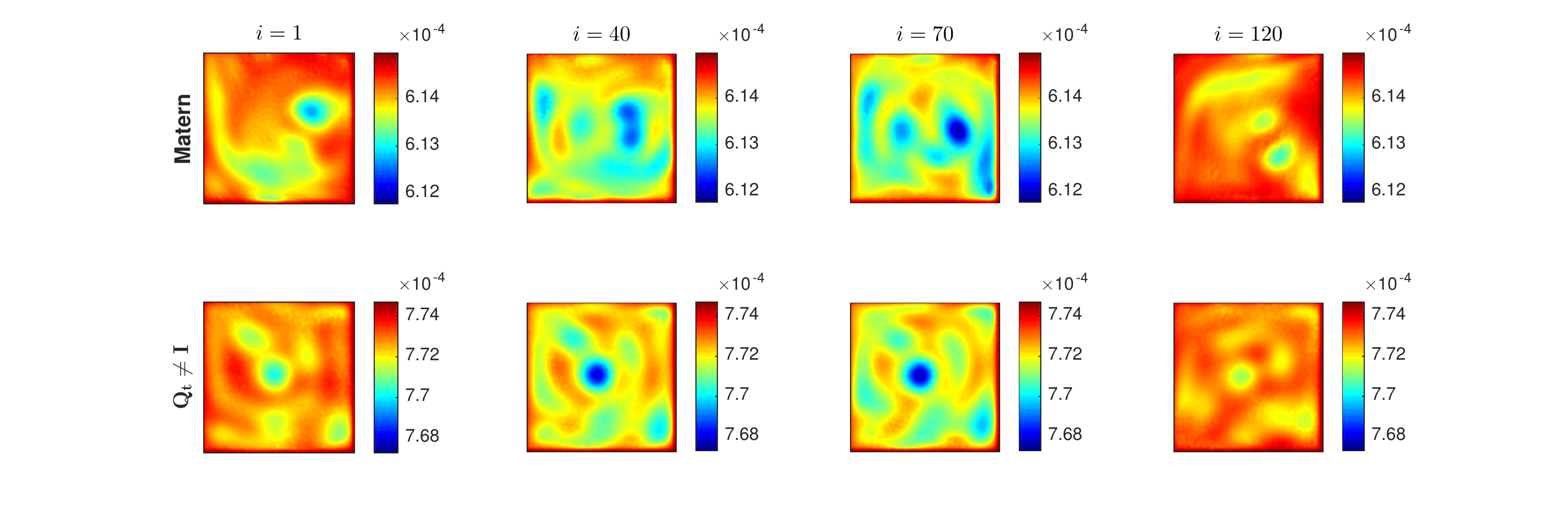}
   \caption{Variance estimates for dynamic PAT reconstructions obtained using the gen-GK bidiagonalization with WGCV after $10$ iterations for various time points ($i=1, 40, 70, 120$).}\label{fig:PAT_var}
 \end{figure}

Next, we show that variance estimates (i.e., approximations to diagonals of the posterior covariance matrix) can be obtained {with minimal additional costs (here, in 15 seconds)}.  In Figure~\ref{fig:PAT_var} we provide results for Mat\'{e}rn and $\bfQ_t \neq \bfI$ (corresponding to the MAP estimates in the 2nd and 4th rows of Figure~\ref{fig:PAT}), where we note that both approaches provide overall variances on the order of $10^{-4}$.  We observe that solutions corresponding to earlier and later time points (e.g., $i=1$ and $i=120$) contain higher variances (i.e., greater uncertainty), with smaller variances in the center regions of the images, especially for $\bfQ_t \neq \bfI.$  Variance images for $\bfQ_t = \bfI$ were essentially constant with mean value $0.0257$ and standard deviation $2.3\times 10^{-6}$ and thus are omitted.

In summary, we have shown that the gen-GK bidiagonalization can be used for the efficient computation of MAP estimates and variance estimates for dynamic PAT problems where the underlying object is changing slowly relative to the rate of image acquisition.  Various choices for the prior covariance matrices could be included in this framework.

 \subsection{Space-time image deblurring}
 \label{sub:deblurring}
\begin{figure}[ht]
   \centering
    \begin{tabular}{cc}
      \includegraphics[width=.4\textwidth]{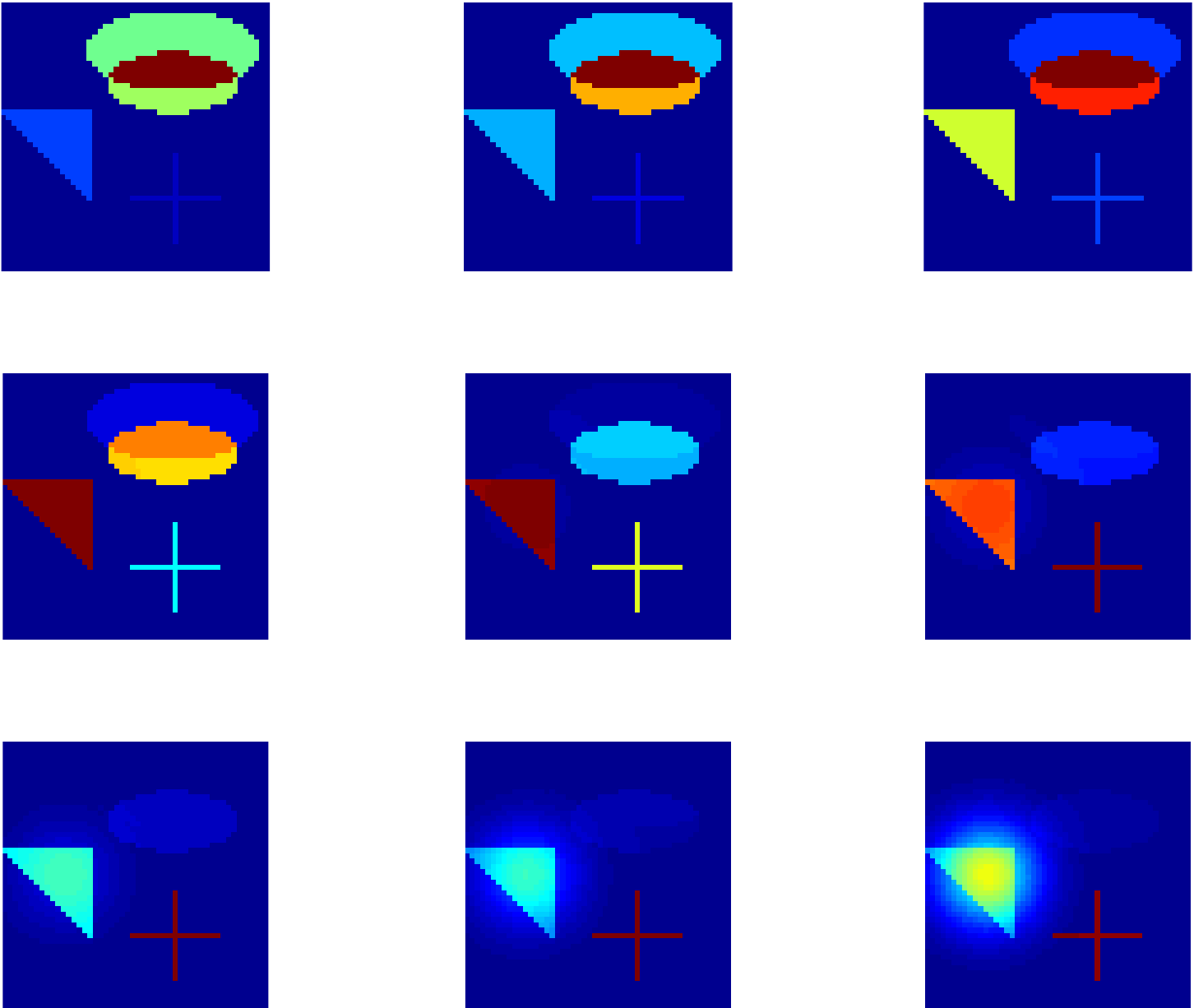} \qquad  & \qquad \includegraphics[width=.4\textwidth]{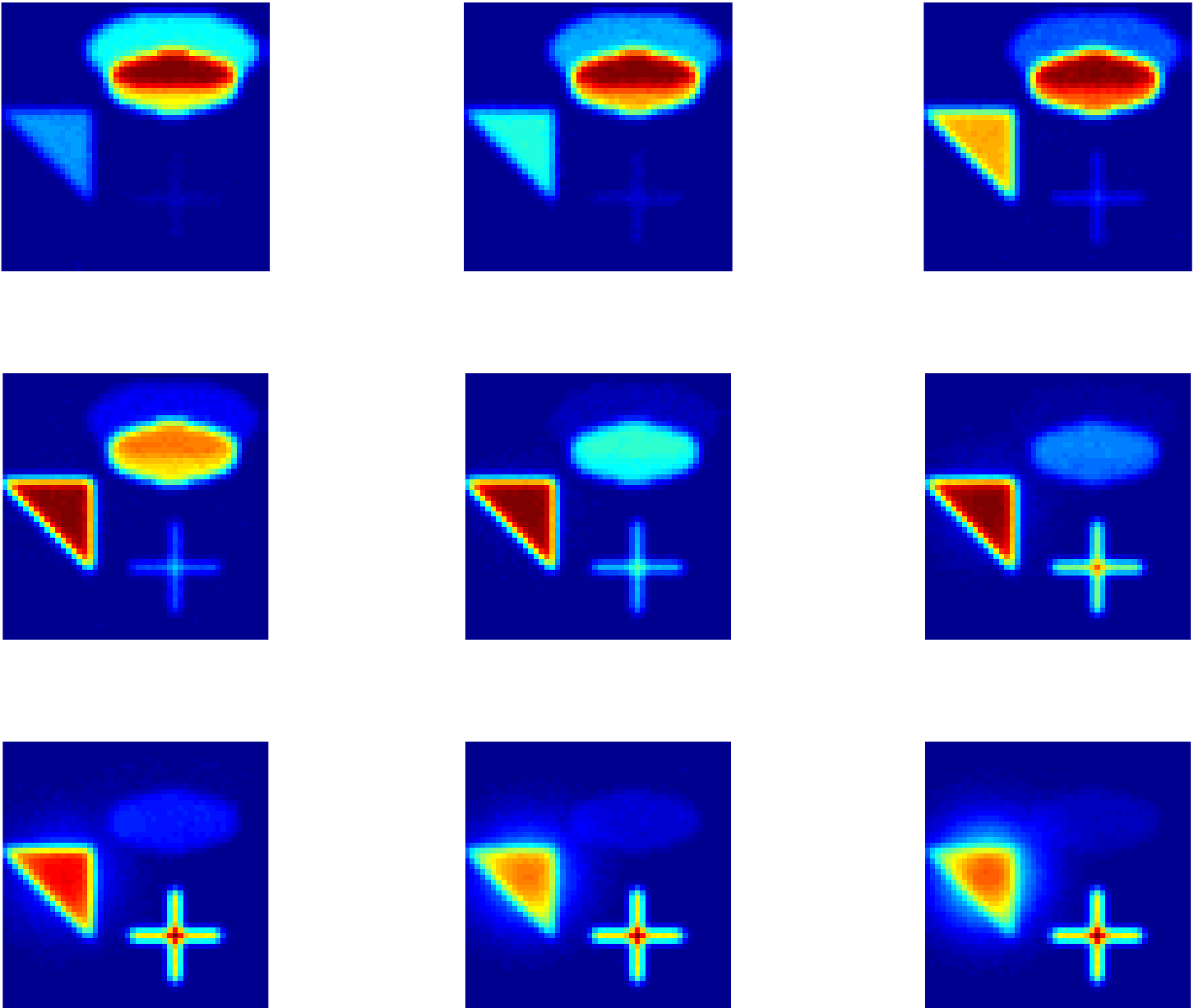}\\
      (a)True images  & (b) Observed, blurred images
    \end{tabular}
   \caption{Dynamic image deblurring problem. }\label{fig:DB_problem}
 \end{figure}

 In dynamic image deblurring, the goal is to reconstruct a sequence of images from a sequence of blurred and noisy images.  We consider a simulated problem where $9$ true images of size $50 \times 50$ are shown in Figure~\ref{fig:DB_problem}(a) and the corresponding observed images are shown in Figure~\ref{fig:DB_problem}(b).  The blur matrix was taken to be
 $\bfA = \bfA_t \kron \bfA_s$ where $\bfA_s$ represents a 2D Gaussian point spread function with spread parameter $\sigma =.07$ and bandwidth $3$ and $\bfA_t$ represents a 1D Gaussian blur with spread parameter $\sigma =1$ and bandwidth $3$.  The noise level is set to be $0.02$ such that $\bfR = 0.0437^2 \bfI.$  The problem set-up is a modification of the `blur' example from the Regularization Tools toolbox \cite{hansen1994regularization}.

\begin{figure}[ht]
  \begin{center}
    \includegraphics[width=\textwidth]{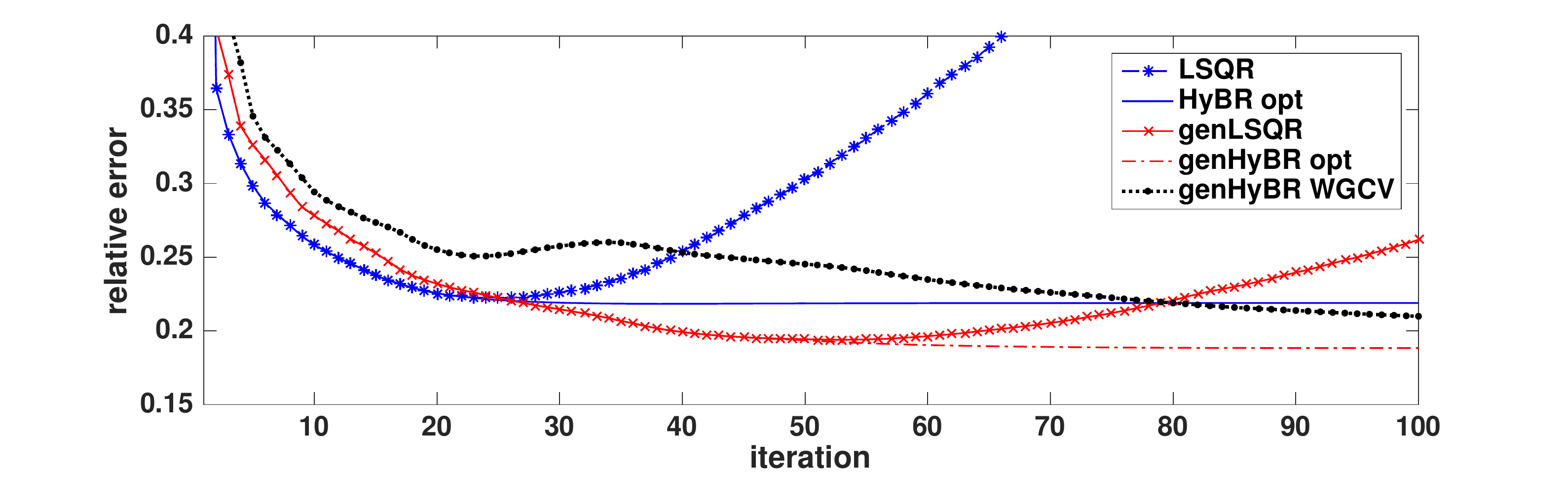}
   \end{center}
   \caption{Relative reconstruction errors for the dynamic image deblurring problem.  LSQR and HyBR correspond to $\bfQ = \bfI$, and genLSQR and genHyBR include $\bfQ=\bfQ_t \kron \bfQ_s$.}\label{fig:DB_err}
 \end{figure}

 We compare LSQR, HyBR opt, genLSQR and genHyBR with two different regularization parameter selection techniques: optimal regularization parameter and the WGCV parameter. Here genLSQR means that $\lambda=0$, and LSQR means that $\lambda=0$ and $\bfQ = \bfI$.
 For genLSQR and genHyBR, we used $\bfQ = \bfQ_t \kron \bfQ_s$ where $\bfQ_t$ and $\bfQ_s$ correspond to $C_T(\cdot)=C_{1.5, .3}(\cdot)$ and $C_S(\cdot)=C_{.5, .007}(\cdot)$ respectively. 
 Relative errors per iteration provided in Figure~\ref{fig:DB_err} reveal similar behavior as that described in~\cite{chungsaibaba2017}.  In particular, LSQR and genLSQR are plagued by semiconvergence (i.e., the ``U''-shaped error curve that results from noise contamination during inversion), which can be avoided in the hybrid variants with the selection of the optimal regularization parameter.  WGCV is able to provide a fairly good regularization parameter, but the process terminated at iteration $24$ due to a flat GCV curve.

 Since $\bfA$ also has a  Kronecker product structure, the decoupled approach applies here. We computed MAP approximations using the decoupled approach, where genHyBR was used to solve each subproblem~\eqref{e_decoupled}.  We denote `decoupled $\lambda_{\rm wgcv}$' to be the solution using the decoupled approach with a fixed regularization parameter $\lambda_{\rm wgcv}$, and `decoupled $\lambda_{\rm wgcv}^{(i)}$' refers to using a different regularization parameter for each subproblem.  WGCV-selected regularization parameters $\lambda_{\rm wgcv}^{(i)}$ and corresponding stopping iterations $k_{\rm stop}^{(i)}$ are provided in Table~\ref{tab:regpar}, along with regularization parameters $\lambda_{\rm wgcv}$ and $\lambda_{\rm opt}$.  We remark that the regularization parameters $\lambda_{\rm wgcv}^{(i)}$ in decoupled approach decrease with increasing index $i$; this can be attributed to the scaling factor from the singular values and the changing right hand sides.
 \begin{table}[ht]
  \centering
   \begin{tabular}{|c|c|c|c|c|c|c|c|c|c|}
     \hline
   $i$  & 1 & 2 & 3 & 4 & 5 &  6 & 7 & 8 & 9\\ \hline
   $\lambda_{\rm wgcv}^{(i)}$ &
   5.4689 &
   5.1159 &
   4.3018 &
   4.0963 &
   2.3008 &
   4.1293 &
   1.5994 &
   0.6247 &
   0.3141 \\ \hline
    $k_{\rm stop}^{(i)}$ &
    61 &
    82 &
    55 &
    24 &
    17 &
     5 &
     5 &
     4 &
     5
     \\ \hline
 $\lambda_{\rm wgcv}$ &   \multicolumn{9}{c| }{8.5838 ($k_{\rm stop} = 22$)}\\ \hline
 $\lambda_{\rm opt}$ &   \multicolumn{9}{c|}{2.2812}   \\ \hline
   \end{tabular}
   \caption{Regularization parameters and stopping iteration in decoupled $\lambda_{\rm wgcv}^{(i)}$.  WGCV regularization parameter (along with stopping iteration for the simultaneous approach) and the optimal regularization parameter are provided for comparison.}
   \label{tab:regpar}

 \end{table}

\begin{figure}[bthp]
  \centering
    \includegraphics[width=.8\textwidth]{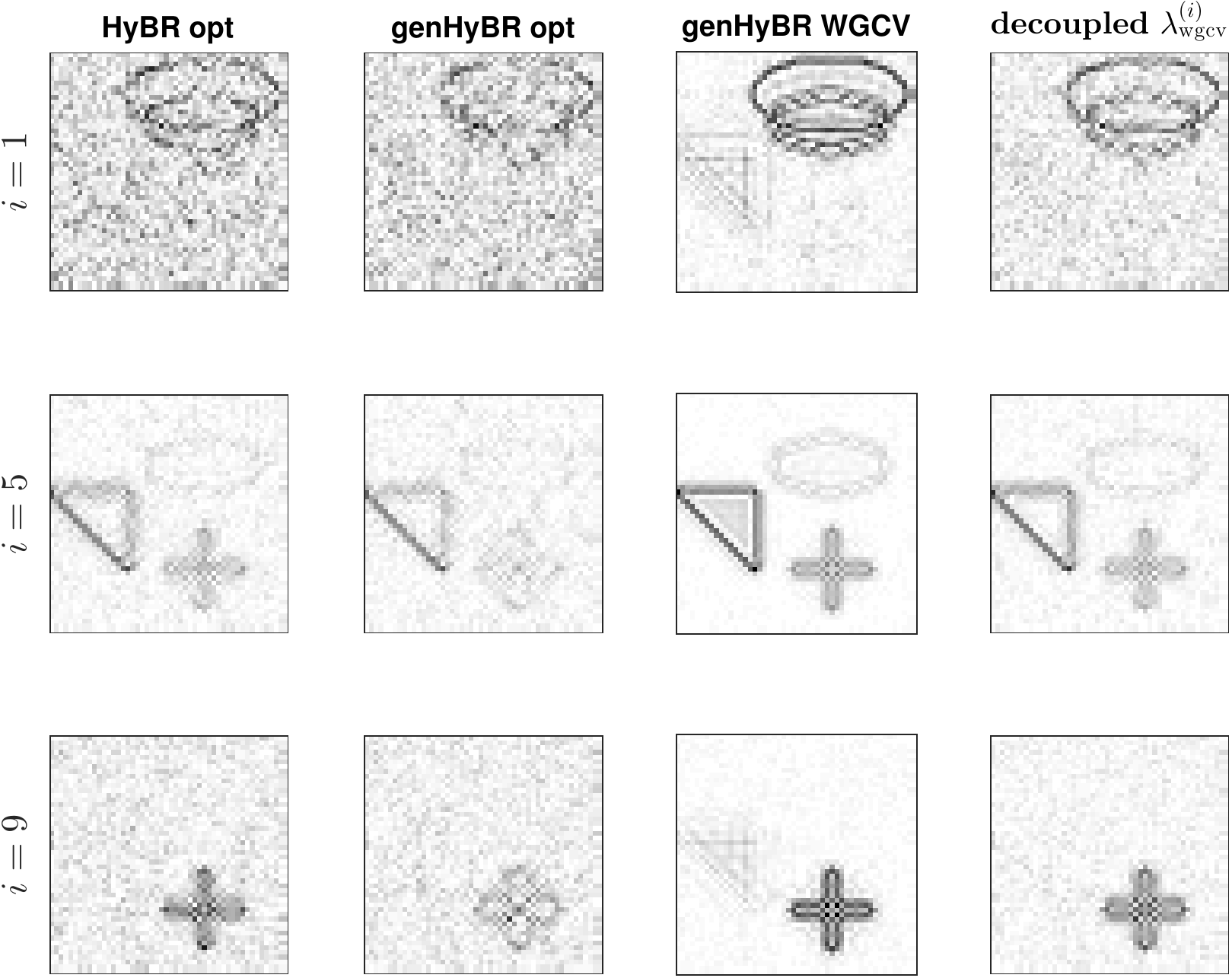}
   \caption{Absolute error images (where white corresponds to a low absolute error and black corresponds to higher error) for slices 1, 5, and 9.}\label{fig:DB_errimages}
 \end{figure}

 The relative reconstruction error for the decoupled $\lambda_{\rm wgcv}$ reconstruction was $0.2461$, which is slightly smaller than that of (simultaneous) genHyBR WGCV which was $0.2507$ at termination.  Furthermore, we observed that allowing different regularization parameters can result in smaller reconstruction error. The relative error for the decoupled $\lambda_{\rm wgcv}^{(i)}$ reconstruction was $0.2198$.  Also, allowing a different Krylov subspace for each reconstruction can be beneficial in reducing ``ghosting'' errors from neighboring images, as evident in the absolute error images provided in Figure~\ref{fig:DB_errimages}, where the same color axis is used per $i$ and black correspond to larger errors.  HyBR opt and genHyBR opt are only provided for reference, since they require the optimal regularization parameter that is not available in practice.

\begin{figure}[bt!]
   \centering
    \includegraphics[width=\textwidth]{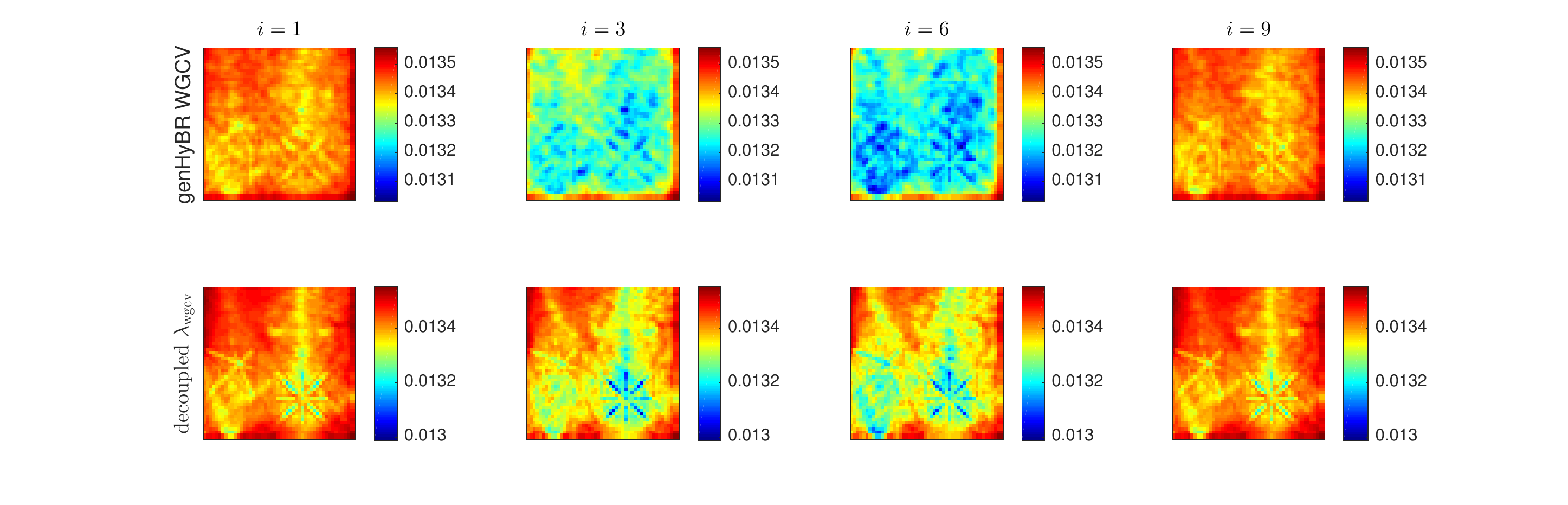}
   \caption{Variance estimates for genHyBR WGCV and decoupled $\lambda_{\rm wgcv}$ for slices 1, 3, 6, and 9. The top row refers to the simultaneous approach, whereas the bottom row corresponds to the decoupled approach.  }\label{fig:DB_varest}
 \end{figure}

 Variance estimates for genHyBR WGCV are provided in the first row of Figure~\ref{fig:DB_varest}, and variance estimates in the second row illustrate that the decoupled variance estimate approach with fixed regularization parameter (as explained in Section~\ref{sub:decvar}) can provide a good approximation.
 We note that decoupled $\lambda_{\rm wgcv}^{(i)}$ does not directly fit our framework; however, a modification of $\bfQ_t^{-1}$ may be used to incorporate the changing regularization parameters.
 In summary, the decoupled approach can be used for both MAP and variance estimation if $\bfA$, $\bfQ,$ and $\bfR$ are all Kronecker products.

As a final remark, in Section~\ref{sub:decoupled} we assume that the regularization parameter is fixed; however, in this section, we present results for which the regularization parameter is allowed to be different. This certainly has benefits as demonstrated; however, its statistical meaning is not fully clear and is worth exploring in future work.

\subsection{Passive seismic tomography (PST)}

Recent advances in PST have enabled the monitoring of mining-induced stress redistribution in coal and hardrock mines \cite{luxbacher2008three,ma2016imaging,westman2012passive}.
The basic goal of microseismic tomography is to image subsurface properties by using the many low-magnitude seismic events (e.g., microearthquakes) that are recorded by a microseismic monitoring system in a deep mine.  Using time-lapse PST tomogram reconstructions, we can better understand the stress redistribution within the rock mass so that trends preceding and following significant seismic events can be analyzed.  However, obtaining these 3D spatial reconstructions in real-time is a computationally challenging task that consists of solving a sequence of very large, often nonlinear, inverse problems.

In this work, we consider a simplified, linear PST problem in a dynamic framework and investigate gen-GK methods for computing reconstructions. The basic formulation of the problem is the same as~\eqref{eq:invproblem} where $\bfs$ is a discretization of the velocity model, $\bfd$ contains the observed travel times or recorded sinogram, and $\bfA$ simulates a ray trace operation.
We consider the situation in which  measurements are taken in periodic time intervals, and the goal is to generate velocity models over time, from which the changing conditions within the rock mass, inferred to be changing due to induced seismicity, can be obtained.

\begin{figure}[ht]
	\centering
\includegraphics[width=\textwidth]{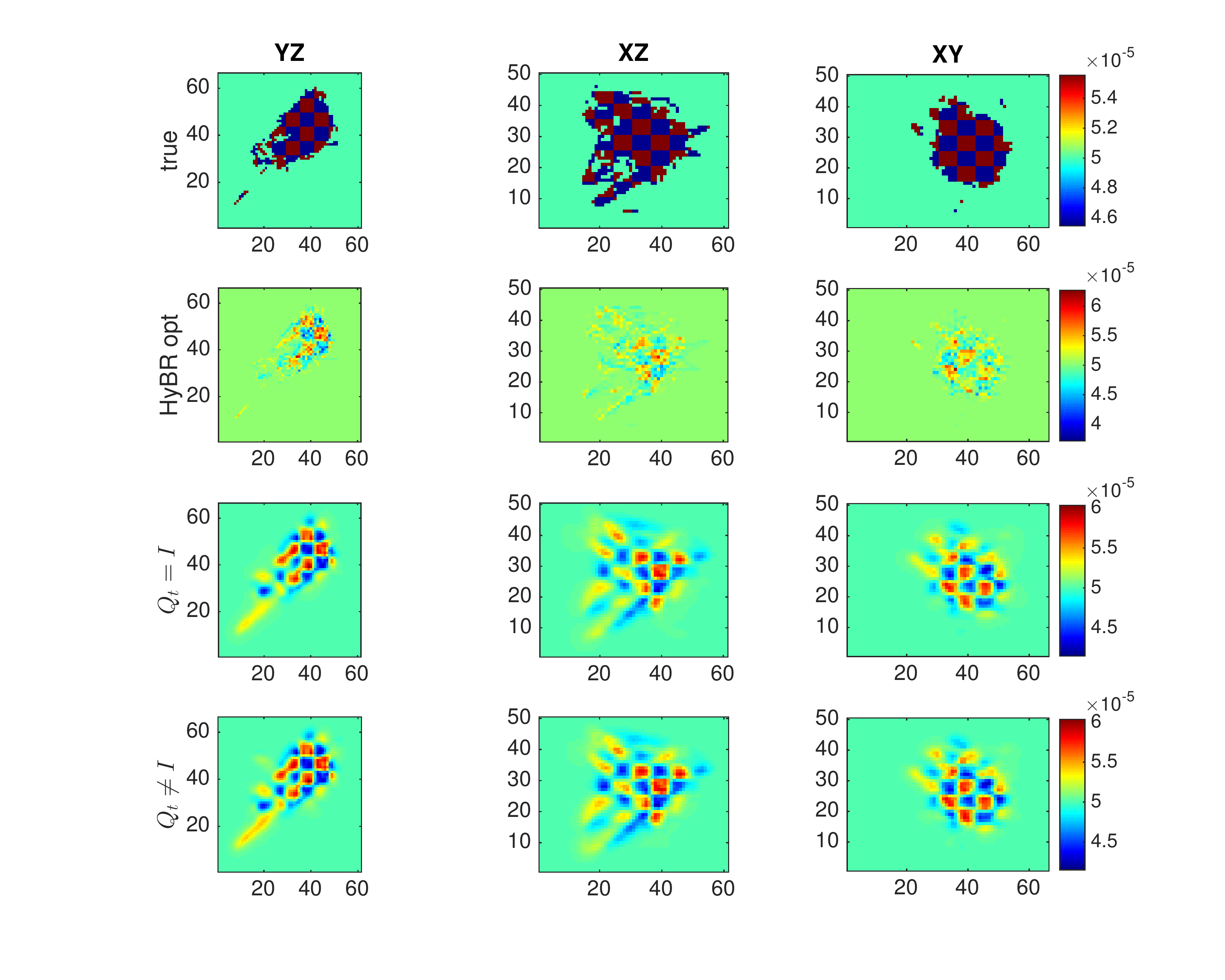}
\caption{Slices from the $4$th volume of the simulated PST example.  The top row contains slices of the true checkerboard volume, the second row contains the HyBR reconstruction (with $\bfQ_t = \bfI$ and $\bfQ_s = \bfI$), and the bottom two rows contain genHyBR reconstructions for $\bfQ_t=\bfI$ and $\bfQ_t \neq \bfI$ respectively. In these images, only the regions with ray coverage have been highlighted.}
\label{fig:checker_recon_mu}
\end{figure}

\paragraph{PST simulated data.}
As is commonly done in practice, we begin with a simulated problem where the goal is to reconstruct a ``checkerboard'' image \cite{westman2012passive}.  We create eight checkerboard volumes of size $50 \times 66 \times 61$ voxels to represent true images, where the values of the checkerboard structure are generated to be reciprocals of $20,000\pm 10\%$ (i.e., values are $4.545\times 10^{-5}$ and $5.555\times 10^{-5}$) in the region of the volume which is seismically ``observable.'' Cross-sections from the $4$th such generated structures are provided in the top row of Figure~\ref{fig:checker_recon_mu}.
Then we used straight-path ray trace matrices $\bfA_i$ for $i=1,\ldots,8$ from real mine data to generate sinograms.  Each matrix corresponds to seismic events that occurred and were detected in a given time period; see Table~\ref{tab:PSTinfo}.  Observed data were constructed using~\eqref{eq:invproblem} where $\bfvarepsilon$ represented Gaussian noise with $\bfR=.0015^2\, \bfI.$  In summary, the number of unknowns for this problem is $50*66*61*8=1,610,400$ and the total number of observations is $m=191,856$ where the number of observations per time point is provided in Table~\ref{tab:PSTinfo}.
\begin{table}[ht]
\begin{center}
	\caption{Dates and number of rays for 8 time periods studied}
	\label{tab:PSTinfo}
	\begin{tabular}{|c|c|c|}\hline
		Begin Date & End Date & Number of rays\\\hline
		Jan 14 & Jan. 24 & $m_1=13604$\\\hline
		Jan 25 & Jan 31 & $m_2=\,\,8512$\\\hline
		Feb 1 & Feb 17 & $ m_3=\,\,9255$\\\hline
		Feb 8 & Feb 14 & $m_4=27121$\\\hline
	 Feb 15 & Feb 21 & $m_5=41282$\\\hline
	 Feb 22 & Feb 28 & $m_6=13755$\\\hline
	 Mar 1 & Mar 17& $m_7=11774$\\\hline
	 Mar 8 & Mar  16 & $m_8=66553$\\ \hline
	\end{tabular}
\end{center}
\end{table}

\begin{figure}[ht]
	\centering
\includegraphics[width=\textwidth]{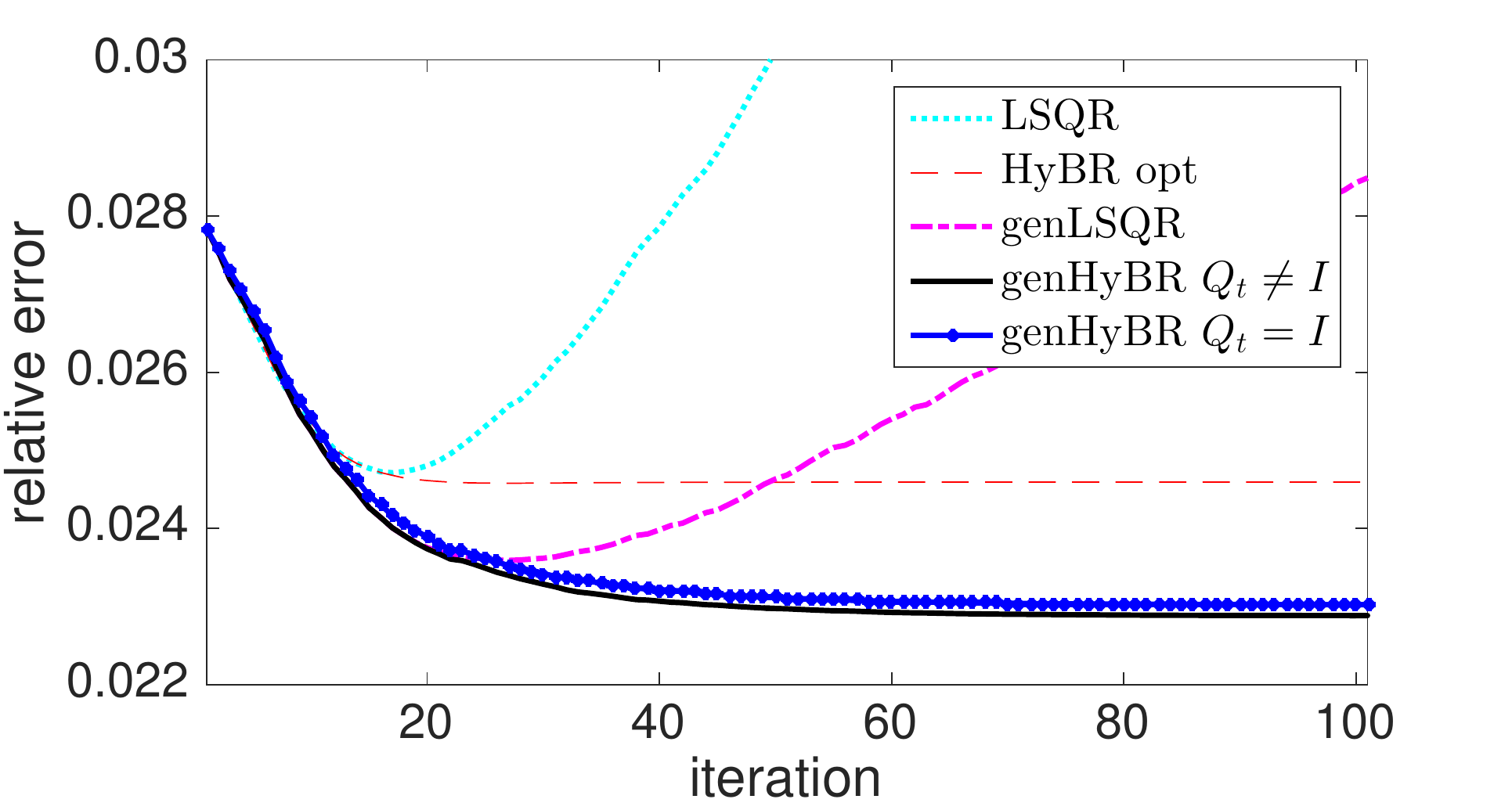}
\caption{Relative reconstruction errors {in the observable regions} per iteration for the PST checkerboard example.}
\label{fig:checker_rel_mu}
\end{figure}

In Figure~\ref{fig:checker_rel_mu}, we provide relative reconstruction errors in the observable regions
per iteration for various methods:
\begin{itemize}
	\item genLSQR corresponds to $C_T(\cdot) = C_{3.5, 0.09}(\cdot)$, $C_S=C_{2.5, .025}(\cdot)$ and $\lambda=0$.
	\item genHyBR $\bfQ_t \neq \bfI$ corresponds to $C_T(\cdot) = C_{3.5, 0.09}(\cdot)$, $C_S=C_{2.5, .025}(\cdot)$.
	\item genHyBR $\bfQ_t = \bfI$ corresponds to $C_S(\cdot)=C_{2.5, .025}(\cdot)$.
\end{itemize}
Here $C_{\nu,\ell}(\cdot)$ corresponds to a covariance matrix determined from a Mat\'ern kernel with parameters $\nu$ and $\ell$ as defined in Section~\ref{sub:modelingtemp}, and we assume, for simplicity,  that $\bfp$ and $t$ are equally spaced on a grid that is normalized to $[0,1].$
Regularization parameter $\lambda_{\rm opt}$ was used for genHyBR.

For LSQR and HyBR (defined in Section~\ref{sub:deblurring}), we follow current practice and take an initial guess $\bfs_0$ to be a constant image with all entries equal to $5\times 10^{-5}$.
For the genHyBR reconstructions, we use a physically-informed prior mean and take $\bfmu$ to be a constant image with all entries equal to $5\times 10^{-5}$.
We found these inclusions to be critical for obtaining physically meaningful results.


Oftentimes in PST, the temporal prior is ignored (e.g., $\bfQ_t = \bfI$) and reconstructions for each time point are done independently (with $\bfQ_s = \bfI$).
Comparing HyBR and genHyBR for $\bfQ_t=\bfI$, we observe that better reconstructions can be obtained by including a spatial prior $\bfQ_s$.  Furthermore, these results show that incorporating a temporal prior may lead to additional improvements.
Cross-sections of the $4$th volume are provided in Figure~\ref{fig:checker_recon_mu} for HyBR, and genHyBR for $\bfQ_t=\bfI$ and $\bfQ_t \neq \bfI$.

\paragraph{PST real data.}  As demonstrated in the simulated problem, the small ray path coverage makes dynamic PST a highly ill-posed problem.  Next, we consider the real field data measurements and reconstruct a time-lapse of eight volumes using genHyBR.  {T}he true volumes are unknown; here we only show isosurfaces and comparisons to currently used algorithms. Along with expert field knowledge, this information can aid in evaluating the reconstructions and the potential for future improvements.

We present three reconstructions using genHyBR.  The first approach essentially mimics what is done in practice, which is to compute reconstructions independent of time.  Here we used genHyBR for each time period with $\bfQ_s$ corresponding to Mat\'ern covariance function $C_S(\cdot)=C_{10.5, .006}(\cdot)$.  In the time independent approach, different regularization parameters and stopping iterations were selected for each reconstruction.  Then we used simultaneous genHyBR with $\bfQ_t = \bfI$ and $\bfQ_t$ corresponding to $C_T(\cdot)=C_{.3, .3}(\cdot)$.  In both cases, we used $\bfQ_s$ as defined above.  For all of these experiments, WGCV was used to select the regularization parameter and automatic stopping criteria was used as described in \cite{ChNaOLe08} with a maximum of 10 iterations.  Obtaining one dynamic PST reconstruction after 10 simultaneous genHyBR iterations took approximately $11$ seconds.

\begin{figure}[ht]
	\centering
	 \includegraphics[width=\textwidth]{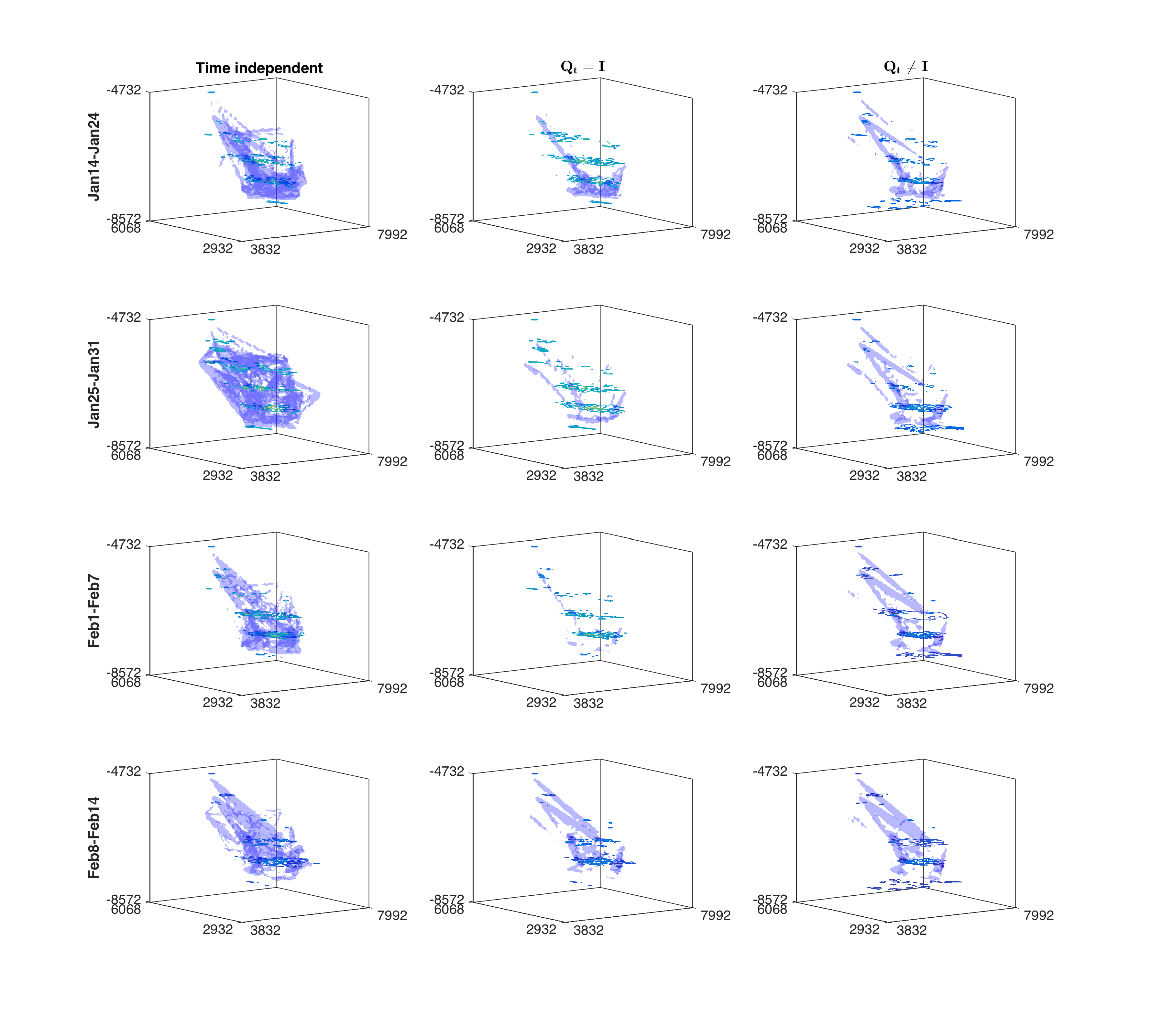}
\caption{Isosurfaces and contours for real PST data at various time points.  The first column corresponds to the time independent genHyBR approach, and the second and third columns correspond to the simultaneous genHyBR approach with $\bfQ_t = \bfI$ and $\bfQ_t \neq \bfI$ respectively.}
\label{fig:recon_real1}
\end{figure}

\begin{figure}[ht]
	\centering
	 \includegraphics[width=\textwidth]{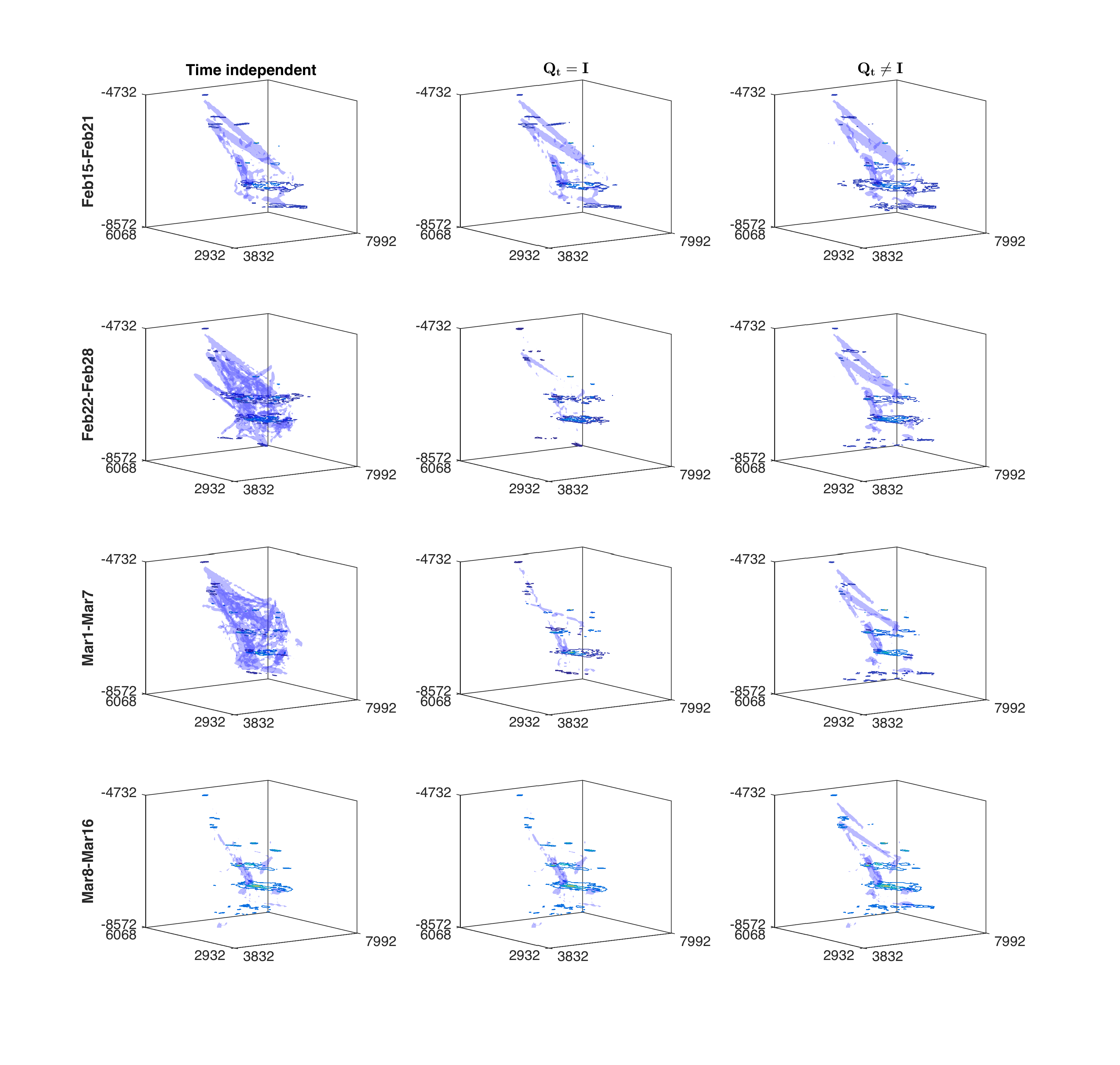}
\caption{Isosurfaces and contours for real PST data at various time points.}
\label{fig:recon_real2}
\end{figure}

High velocity iso-values (corresponding to a value of $20,025$ or a slowness of $4.9938\times 10^{-5}$) and contours for different time intervals are shown in Figures~\ref{fig:recon_real1} and~\ref{fig:recon_real2}.  Notice that including the temporal prior can result in better reconstructions, especially for time intervals with very few observations (e.g., Jan 25--31 and Feb 1--17).  In addition, genHyBR with $\bfQ_t \neq \bfI$ gives more detailed information and locality of high stresses than the time independent reconstructions, where the isosurfaces are more dispersed.

\begin{figure}[ht]
	\centering
	 \includegraphics[width=\textwidth]{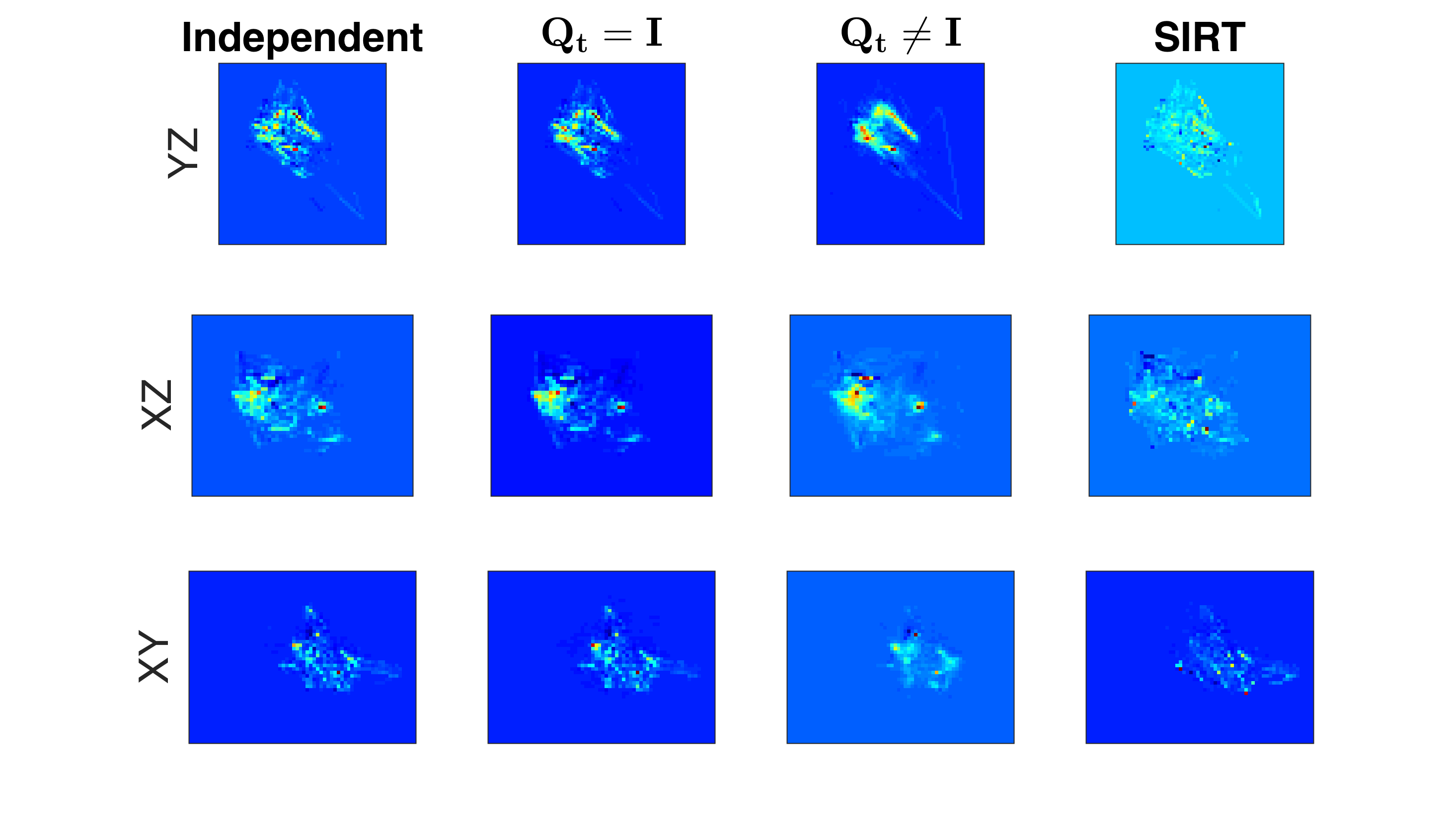}
\caption{Cross-sections of reconstructions for real PST data.}
\label{fig:recon_sirt}
\end{figure}

 We compare our reconstructions to a typical field reconstruction using the simultaneous iterative reconstruction technique (SIRT) for the time period February 22--28. In Figure~\ref{fig:recon_sirt} we display cross-sections of genHyBR reconstructions, along with the SIRT reconstruction.  These correspond to the $24$th, $40$th, and $27$th slices in the $x$, $y$, and $z$ dimensions respectively.  We remark that the genHyBR reconstructions all provide smoother, more localized reconstructions of high-velocity zones.  It is worth mentioning that these reconstructions have assumed that the forward model has  ``straight-ray'' paths and a typical approach in mining would be to use this reconstruction as an initial guess for {obtaining} reconstructions with a more sophisticated, nonlinear ``curved-ray'' forward model.  This is another topic of future work.

\section{Conclusions}
\label{sec:conclusions}

We consider the problem of dynamic inverse problems using the Bayesian framework, and we developed efficient, iterative, matrix-free methods based on the gen-GK bidiagonalization.  A wide range of priors can be incorporated in our framework.  We focused on priors that are modeled as Gaussian random fields with special attention to space-time covariance kernels for which MVPs can be computed efficiently.
We first focus on computing the MAP estimate. In the simultaneous approach, a solution for the entire unknown in space-time is solved in an ``all-at-once'' manner. When the observation operator also has Kronecker product structure, a series of variable transformations enables the problem to decouple in time. Both the simultaneous and decoupled approaches leverage the efficient iterative solvers developed in our previous work~\cite{chungsaibaba2017}, and has the added benefit that the simultaneous approach allows for automatic selection of regularization parameter selection. In addition to the MAP estimate, we describe methods that reuse intermediate information contained in the iterative solvers to estimate the variance of the posterior distribution.  Several examples from image processing, including new applications to PST, demonstrate scalability of our algorithms and illustrate the broad applicability of our work.

\section{Acknowledgements}
We would like to acknowledge Creighton mines for providing the raw data for the {PST} application.  Furthermore, some of this work was conducted as a part of the SAMSI Program on Optimization 2016-2017. This material was based upon work partially supported by the National Science Foundation under Grant DMS-1127914 to the Statistical and Applied Mathematical Sciences Institute.

\bibliography{7references-edit}
\bibliographystyle{plain}

\end{document}